\RequirePackage{tikz}

\documentclass[sn-basic,iicol]{sn-jnl}

\newcommand{\bA}{\boldsymbol{A}}
\newcommand{\bI}{\boldsymbol{I}}
\newcommand{\bL}{\boldsymbol{L}}
\newcommand{\bvarepsilon}{\boldsymbol{\varepsilon}}
\newcommand{\bsigma}{\boldsymbol{\sigma}}
\newcommand{\TP}{^\text{T}}
\newcommand{\eMod}{E}  
\newcommand{\poission}{\nu}  
\newcommand{\poissionZero}{\poission^{(0)}}  
\newcommand{\fc}{f_{\text{c}}}  
\newcommand{\fcInf}{{\fc}_\infty}  
\newcommand{\fcEff}{{\fc}_{,\text{eff}}}  
\newcommand{\body}{\Omega}  
\newcommand{\phaseIndex}{r}  
\newcommand{\matrixIndex}{\text{m}}  
\newcommand{\inclIndex}{\text{i}}  
\newcommand{\bodyPhase}{\body^{(\phaseIndex)}}
\newcommand{\volFrac}{c}
\newcommand{\volFracPhase}{\volFrac^{(\phaseIndex)}}
\newcommand{\volFracMatrix}{\volFrac^{(\matrixIndex)}}
\newcommand{\volFracIncl}{\volFrac^{(\inclIndex)}}
\newcommand{\volFracZero}{\volFrac^{(0)}}
\newcommand{\radius}{R}
\newcommand{\radiusPhase}{\radius^{(\phaseIndex)}}
\newcommand{\bulkMod}{K}
\newcommand{\bulkModPhase}{\bulkMod^{(\phaseIndex)}}

\newcommand{\bulkModZero}{\bulkMod^{(0)}}
\newcommand{\bulkModEff}{\bulkMod_{\text{eff}}}
\newcommand{\shearMod}{G}
\newcommand{\shearModPhase}{\shearMod^{(\phaseIndex)}}

\newcommand{\shearModZero}{\shearMod^{(0)}}
\newcommand{\shearModEff}{\shearMod_{\text{eff}}}
\newcommand{\matStiff}{\bL}
\newcommand{\matStiffZero}{\matStiff^{(0)}}
\newcommand{\matStiffPhase}{\matStiff^{(\phaseIndex)}}
\newcommand{\matStiffEff}{\matStiff_{\text{eff}}}
\newcommand{\orthProjV}{\bI_{\text{V}}}
\newcommand{\orthProjD}{\bI_{\text{D}}}
\newcommand{\strain}{\bvarepsilon}
\newcommand{\strainPhase}{\strain^{(\phaseIndex)}}
\newcommand{\strainZero}{\strain^{(0)}}
\newcommand{\stress}{\bsigma}
\newcommand{\stressZero}{\stress^{(0)}}
\newcommand{\stressMatrix}{\stress^{(\matrixIndex)}}
\newcommand{\stressD}{\stress_{\text{D}}}
\newcommand{\stressTest}{\stress^{\text{test}}}
\newcommand{\concentration}{\bA}
\newcommand{\concentrationPhase}{\concentration^{(\phaseIndex)}}
\newcommand{\concentrationZero}{\concentration^{(0)}}
\newcommand{\dilConcentration}{\concentration_{\text{dil}}}
\newcommand{\dilConcentrationPhase}{\dilConcentration^{(\phaseIndex)}}

\newcommand{\dilConcentrationVPhase}{A^{(\phaseIndex)}_{\text{dil,V}}}
\newcommand{\dilConcentrationDPhase}{A^{(\phaseIndex)}_{\text{dil,D}}}

\newcommand{\auxAlphaZero}{\alpha^{(0)}}
\newcommand{\auxBetaZero}{\beta^{(0)}}
\newcommand{\Jtwo}{J_2}
\newcommand{\JtwoTest}{\Jtwo^{\text{test}}}

\newcommand{\JtwoMatrix}{\Jtwo^{(\matrixIndex)}}
\newcommand{\force}{f} 
\newcommand{\forceTest}{\force^{\text{test}}} 
\newcommand{\thermCond}{\chi}

\newcommand{\thermCondMatrix}{\thermCond^{(\matrixIndex)}}
\newcommand{\thermCondIncl}{\thermCond^{(\inclIndex)}}
\newcommand{\thermCondHom}{\thermCond_{\text{eff}}}

\newcommand{\concentrationThermCondIncl}{A_{\chi}^{(\inclIndex)}}
\newcommand{\Ivol}{\bI_{\text{V}}}
\newcommand{\Idev}{\bI_{\text{D}}}
\newcommand{\currentn}{n+1}
\newcommand{\lastn}{n}
\newcommand{\DOH}{\alpha} 
\newcommand{\DOHLast}{\DOH^{\lastn}} 
\newcommand{\DOHCurrent}{\DOH^{\currentn}} 
\newcommand{\DOHmax}{\alpha_{\text{max}}} 
\newcommand{\DOHt}{\alpha_{\text{t}}} 
\newcommand{\DOHTwentyEight}{\alpha_{28}} 
\newcommand{\heat}{Q} 
\newcommand{\heatInf}{Q_{\infty}} 
\newcommand{\zeit}{t} 
\newcommand{\temp}{T} 
\newcommand{\tempMax}{T_{\text{max}}} 
\newcommand{\tempLimit}{T_{\text{limit}}} 
\newcommand{\tempCurrent}{\temp^{\currentn}} 
\newcommand{\tempLast}{\temp^{\lastn}} 
\newcommand{\tempRef}{\temp_{\text{ref}}}
\newcommand{\dTdt}{\frac{\partial \temp}{\partial \zeit}}  
\newcommand{\dQdt}{\frac{\partial \heat}{\partial \zeit}}  
\newcommand{\heatCapSpecific}{C} 
\newcommand{\density}{\rho} 
\newcommand{\thermCondEff}{\lambda} 
\newcommand{\dDOHdt}{\frac{\partial \DOH}{\partial \zeit}}  
\newcommand{\affinity}{A} 
\newcommand{\affinityTemp}{\tilde{\affinity}} 
\newcommand{\affinityScale}{a} 
\newcommand{\hydParBone}{B_1}
\newcommand{\hydParBtwo}{B_2}
\newcommand{\hydParEta}{\eta}
\newcommand{\function}{f}

\newcommand{\strengthExp}{a}
\newcommand{\fcTwentyEight}{{\fc}_{28}}

\newcommand{\strengthCExp}{\strengthExp_{\fc}}

\newcommand{\stiffExp}{\strengthExp_{\eMod}}
\newcommand{\eModInf}{\eMod_\infty}  
\newcommand{\eModTwentyEight}{\eMod_{28}}  
\newcommand{\activE}{E_{\text{a}}}
\newcommand{\gasConst}{R}
\newcommand{\wc}{r_{\text{wc}}}
\newcommand{\stressInvariant}{I}
\newcommand{\firstStressInvariant}{\stressInvariant_{1}}
\newcommand{\secondStressInvariant}{\stressInvariant_{2}}
\newcommand{\principalStressTension}{\bsigma_{\text{t}}'}
\newcommand{\beamLength}{l}
\newcommand{\beamDistrLoad}{q}
\newcommand{\beamPointLoad}{F}
\newcommand{\beamMaxMoment}{M_{\text{max}}}

\newcommand{\beamHeight}{h}
\newcommand{\beamWidth}{b}
\newcommand{\beamHeightEff}{\beamHeight_{\text{eff}}}

\newcommand{\beamCover}{c}
\newcommand{\beamSteelDiameter}{d}
\newcommand{\beamSteelDiameterMax}{d_{\text{max}}}
\newcommand{\beamSteelDiameterStirrups}{d_{\text{st}}}
\newcommand{\beamConcreteSF}{\gamma_{\text{c}}}
\newcommand{\beamSteelSF}{\gamma_{\text{s}}}
\newcommand{\beamTimeSF}{\alpha_{\text{cc}}}
\newcommand{\beamfcd}{f_{\text{cd}}}
\newcommand{\beamfsd}{f_{\text{ywd}}}
\newcommand{\beamfs}{f_{\text{yk}}}
\newcommand{\beamLeverMoment}{z}
\newcommand{\beamK}{\mu}
\newcommand{\beamSteelReq}{A_{\text{s,req}}}

\newcommand{\beamSteelMax}{A_{\text{s,max}}}
\newcommand{\beamNSteel}{n_{\text{s}}}
\newcommand{\beamNSteelMax}{n_{\text{s,max}}}

\newcommand{\constraint}{\mathcal{C}}
\newcommand{\beamConstraintFc}{\constraint_{\text{fc}}}
\newcommand{\beamConstraintGeometry}{\constraint_{\text{g}}}
\newcommand{\beamConstraintBeam}{\constraint_{\text{beam}}}
\newcommand{\beamMinSpacing}{s_{\text{min}}}
\newcommand{\FEMConstraintT}{\constraint_{\text{T}}}
\newcommand{\FEMConstraintStress}{\constraint_{\stress}}
\newcommand{\FEMConstraintDP}{\constraint_{\text{DP}}}
\newcommand{\FEMConstraintRK}{\constraint_{\text{RK}}}

\newcommand{\pasteE}{30e9}
\newcommand{\pasteEunit}{Pa}
\newcommand{\pastenu}{0.2}

\newcommand{\pasteC}{870}
\newcommand{\pasteCunit}{\dfrac{J}{kgK}}
\newcommand{\pastekappa}{1.8}
\newcommand{\pastekappaunit}{\dfrac{W}{mK}}
\newcommand{\pasterho}{2400}
\newcommand{\pasterhounit}{\dfrac{kg}{m^3}}
\newcommand{\pastefc}{30e6}
\newcommand{\pastefcunit}{Pa}
\newcommand{\pasteQ}{250000}
\newcommand{\pasteQunit}{\dfrac{J}{kg}}
\newcommand{\aggregatesE}{25e9}

\newcommand{\aggregatesnu}{0.3}

\newcommand{\aggregatesC}{840}

\newcommand{\aggregateskappa}{0.8}

\newcommand{\aggregatesrho}{2600}

\newcommand{\beamExSpanUnit}{cm}

\newcommand{\beamExWidthUnit}{mm}

\newcommand{\beamExYieldStrSteelUnit}{N/mm^2}
\newcommand{\beamExSteelDiaBu}{10}
\newcommand{\beamExSteelDiaBuUnit}{mm}
\newcommand{\beamExCoverMin}{2.5}
\newcommand{\beamExCoverMinUnit}{cm}
\newcommand{\beamExHeightC}{450}
\newcommand{\beamExHeightCUnit}{mm}
\newcommand{\beamExPointLoadC}{50}
\newcommand{\beamExPointLoadCUnit}{kN}
\newcommand{\beamExComprStrConcreteC}{40}
\newcommand{\beamExComprStrConcreteCUnit}{N/mm^2}

\newcommand{\evoExAlphaT}{0.2}
\newcommand{\evoExAlphaTx}{0.8}
\newcommand{\evoExaE}{0.5}
\newcommand{\evoExafc}{0.5}
\newcommand{\evoExE}{50}
\newcommand{\evoExfc}{30}  

\newcommand{\ee}{\end{equation}}
\newcommand{\be}{\begin{equation}}
\newcommand{\ec}{\end{center}}
\newcommand{\bc}{\begin{center}}
\newcommand{\eea}{\end{eqnarray}}
\newcommand{\bea}{\begin{eqnarray}}
\newcommand{\bd}{\begin{description}}
\newcommand{\ed}{\end{description}}
\newcommand{\bi}{\begin{itemize}}
\newcommand{\ei}{\end{itemize}}
\newcommand{\pa}{\partial}

\newcommand{\bx}{\bs{x}}
\newcommand{\bb}{\bs{b}}

\newcommand{\bs}{\boldsymbol}

\newcommand{\bz}{\bs{z}}

\newcommand{\bt}{\bs{\theta}}

\newcommand{\E}{\mathbb{E}}
\newcommand\calC{\mathcal{C}}
\newcommand\calO{\mathcal{O}}
\newcommand\calL{\mathcal{L}}
\newcommand\calD{\mathcal{D}}

\newcommand{\unit}[1]{\ensuremath{\, \mathrm{#1}}}

\newcommand{\refeq}[1]{Eq. (\ref{#1})} 
\newcommand{\refsec}[1]{Section \ref{#1}}
\newcommand{\reffig}[1]{Fig. \ref{#1}}

\newcommand{\inputtemperaturelimit}{60}

\newcommand{\inputconcretecover}{2.5}
\newcommand{\inputconcretecoverunit}{cm}

\usepackage{subfigure}
\usepackage{amsmath,bm}
\usepackage{amssymb}
\usepackage{mdframed}
\usepackage{array} 
\usepackage{tikz}
\usetikzlibrary{graphs,graphs.standard}
\usepackage[labelsep=space]{caption}
\usepackage{pdfcomment}
\usepackage{changes}
\usepackage{stanli}

\usepackage{graphicx}%
\usepackage{multirow}%
\usepackage{amsmath,amssymb,amsfonts}%
\usepackage{amsthm}%
\usepackage{mathrsfs}%
\usepackage[title]{appendix}%
\usepackage{xcolor}%
\usepackage{textcomp}%
\usepackage{manyfoot}%
\usepackage{booktabs}%
\usepackage{algorithm}%
\usepackage{algorithmicx}%
\usepackage{algpseudocode}%
\usepackage{listings}%

\usepackage{gensymb}
\usepackage{comment}

\begin{document}
\title[From concrete mixture to structural design - a holistic optimization procedure]{From concrete mixture to structural design - a holistic optimization procedure in the presence of uncertainties}


\author*[1]{\fnm{Atul} \sur{Agrawal}}\email{atul.agrawal@tum.de}
\equalcont{These authors contributed equally to this work.}

\author[2]{\fnm{Erik} \sur{Tamsen}}
\equalcont{These authors contributed equally to this work.}

\author[1]{\fnm{Phaedon-Stelios} \sur{Koutsourelakis}}\email{p.s.koutsourelakis@tum.de}

\author[2]{\fnm{Jörg F.} \sur{Unger}}\email{joerg.unger@bam.de}

\affil*[1]{\orgdiv{Data-driven Materials Modeling}, \orgname{Technische Universität München}, \orgaddress{\street{Boltzmannstraße 15}, \city{Garching}, \postcode{85748}, \country{Germany}}}

\affil[2]{\orgdiv{Modeling and Simulation}, \orgname{Bundesanstalt für Materialforschung und -prüfung}, \orgaddress{\street{Unter den Eichen 87}, \city{Berlin}, \postcode{12205}, \country{Germany}}}

\abstract{
Designing civil structures such as bridges, dams or buildings is a complex task requiring many synergies from several experts. Each is responsible for different parts of the process.
This is often done in a sequential manner, e.g. the structural engineer makes a design under the assumption of certain material properties (e.g. the strength class of the concrete), and then the material engineer optimizes the material with these restrictions.
This paper proposes a holistic optimization procedure, which combines the concrete mixture design and structural simulations in a joint, forward workflow that we ultimately seek to invert. In this manner,  new mixtures beyond standard ranges can be considered. 
Any design effort should account for the presence  of uncertainties which can be aleatoric or epistemic as when data is used to calibrate  physical models or identify models that  fill  missing links in the workflow. 
Inverting the causal relations established poses several challenges especially when these involve physics-based models which most often than not do not provide derivatives/sensitivities or when design constraints are present.
To this end, we advocate Variational Optimization, with proposed extensions and appropriately chosen heuristics to overcome the aforementioned challenges.
The proposed methodology is illustrated using the design of a precast concrete beam with the objective to minimize the global warming potential while satisfying  a number of constraints associated with its  load-bearing capacity after 28days according to the Eurocode,  the demoulding time  as computed by a complex nonlinear Finite Element model, and  the maximum temperature during the hydration. 
}


\maketitle

\section{Introduction}\label{sec:introduction}
Precast concrete elements play a critical role in achieving efficient, low cost and sustainable structures.
The controlled manufacturing environment allows for higher quality products and enables the mass production of such elements.
In the standard design approach, engineers or architects select the geometry of a structure, estimate the loads, choose mechanical properties, and design the element accordingly. 
If the results are not satisfactory, the required mechanical properties are iteratively adjusted, aiming to improve the design.
This approach is adequate  when the choice of mixtures is limited and the expected concrete properties are well known.
There are various published methods to automate this process and optimize the beam design at this level.
Computer-aided beam design optimization dates back at least 50 years, e.g. \citep{Haung1967}.

Generally, the objective is to reduce costs, with the design variables being the beam geometry, the amount and location of the reinforcement and the compressive strength of the concrete, \citep{Chakrabarty_1992, Coello_1997, Pierott_2021, Shobeiri_2023}.
Most publications focus on analytical functions based on well-known, empirical rules of thumb.
In recent years, the use of alternative binders in the concrete mix design has increased, mainly to reduce the environmental impact and cost of concrete but also to improve and modify specific properties.
This is a challenge as the concrete mix is no longer a constant and is itself subject to an optimization.
Known heuristics might no longer apply to the new materials and old design approaches might fail to produce optimal results.
In addition, it is not desirable to choose from a predetermined set of possible mixes, as this would either lead to an overwhelming number of required experiments or a limiting subset of the possible design space.

In the existing literature on  the optimization of the concrete mix design  \citep{Lisienkova_2021, Kondapally_2022}, the  objective is to either improve mechanical properties like durability within constraints, or to minimize e.g. the amount of concrete while keeping other properties above a threshold.
A first step to address these limitations is incorporating the compressive strength during optimization in the beam design phase.
Higher compressive strength usually correlates with a larger amount of cement and, therefore higher cost as well as a higher Global Warming Potential (GWP).
This approach has shown promising results in achieving improved structural efficiency while considering environmental impact \citep{dos_Santos_2023}.
To be able to find a part specific optimum, individual data of the manufacturer and specific mix options must be integrated.
Therefore, there is still a need for a comprehensive optimization procedure that can seamlessly integrate concrete mix design and structural simulations, ensuring structurally sound and buildable elements with minimized environmental impact for part specific data.

When designing elements subjected to various requirements, both on the material and structural level, including workability of the fresh concrete, durability of the structure, maximum acceptable temperature, minimal cost and Global Warming Potential (GWP), the optimal solution is not apparent and will change depending on each individual project.

In this paper, we present a holistic optimization procedure that combines physics-based models and  experimental data in order to enable the optimization of  the concrete mix design in the presence of uncertainty, with an objective to minimize the global warming potential. 
In particular, we employ structural simulations as constraints to ensure structural integrity, limit the maximum temperature and ensure an adequate time of demolding.
 
By integrating the concrete mixture optimization and structural design processes, engineers can tailor the concrete properties to meet specific requirements of the customer and manufacturer.
This approach opens up possibilities for performance prediction and optimization for new mixtures that fall outside the standard range of existing ones.
To the best of our knowledge there are no published works that combine the material and structural level in one flexible optimization framework.
In addition to changing the order of the design steps, the proposed framework allows to directly integrate experimental data and propagate the identified uncertainties.
This allows a straightforward integration of new data and quantification of uncertainties regarding the predictions.
The proposed framework consists of three main parts.
First, an automated and reproducible probabilistic machine learning based parameter identification method to calibrate the models by using experimental data.
Second, a black-box optimization method for non-differentiable functions, including constraints.
Third, a flexible workflow combining the models and functions required for the respective problem.


To carry out black-box optimization, we advocate the use of  Variational Optimization \citep{bird_stochastic_2018,staines2013optimization} which uses stochastic gradient estimators for black-box functions. We utilize this with appropriate enhancements in order  to account for the stochastic, non-linear constraints. Our choice is motivated by three  challenges present in the workflow describing the physical process. Firstly, the  availability of only black-box evaluations of the physical workflow. 
In many real world cases involving physics based solvers/simulators in the optimization process, one resorts to gradient-free optimization \citep{more2009benchmarking,snoek2012practical}.
However, the gradient-free methods perform poorly on high-dimensional parametric spaces \citep{more2009benchmarking}. Also, it requires more functional evaluations to reach the optimum as compared to gradient-based methods.
Recently, stochastic gradient estimators \citep{mohamed2020monte} have been used to estimate gradients of black-box functions and, hence, perform gradient-based optimization \citep{louppe_adversarial_2019,shirobokov2020black,ruiz2018learning}. However, they do not account for the constraints. Secondly, the presence of non-linear constraints. Popular gradient-free methods like constrained Bayesian Optimization (cBO) \citep{gardner2014bayesian} and COBYLA \citep{powell1994direct} pose significant challenges when (non-)linear constraints are involved \citep{menhorn2017trust,audet2016blackbox,agrawal_opt}. 
Thirdly, the stochasticity in the workflow, discussed in the following paragraph.

The physical workflow comprising physics-based models to link design variables with the objective and constraints poses an information flow-related challenge. Some links leading to the objective/constraints are not known a priori in the literature, thus hindering the optimization process. We propose a method to learn these missing links, parameterized by an appropriate neural network, with the help of (noisy) experimental data and physical models. The unavoidable noise in the data introduces aleatoric uncertainty, or its incompleteness introduces epistemic uncertainty. To account for the presence of these uncertainties, we advocate the links to be \textit{probabilistic}. 
The learned probabilistic links tackle the information bottleneck, however, it introduces random parameters in the physical workflow, thus necessitating Optimization under uncertainty (OUU) \citep{martins2021engineering,acar2021modeling}. Deterministic inputs can lead to a poor-performing design, which OUU tries to tackle by producing a robust and reliable design that is less sensitive to inherent variability. 
This paradigm of fusing data and physical models to train machine-learning models has been extensively used across engineering and physics in recent years \citep{karpatne2022knowledge,lucor2022simple,koutsourelakis2016big,agrawal2023probabilistic,fleming2018artificial,karniadakis2021physics}, colloquially referred to as Scientific Machine Learning (SciML). In contrast to traditional machine learning areas where big data is generally available, engineering and physical applications generally suffer from a lack of data, further complicated by experimental noise. Scientific Machine Learning has shown promise in addressing this lack of data.

The structure of the rest of the paper is as follows. \refsec{sec:design_approach} describes the proposed design approach, \refsec{sec:models} describes the physical material models and the applied assumptions. \refsec{sec:exp_data} presents the details of the experimental data. \refsec{sec:calibration} provides an overview of the aforementioned probabilistic links and the optimization procedure. \refsec{sec:model_learning} talks about the methodology employed to learn the probabilistic links based on the experimental data and the physical models. Then \refsec{sec:optimization} describe the details of the proposed black-box optimization algorithm. In \refsec{sec:numericalexperiments}, we showcase and discuss the results of the numerical experiments combining all the parts, the experimental data, the physical models, the identification of the probabilistic links, and the optimization framework. Finally, in \refsec{sec:conclusion}, we summarize our findings and discuss possible extensions.

\subsection{Demonstration  problem}\label{sec:example_problem}
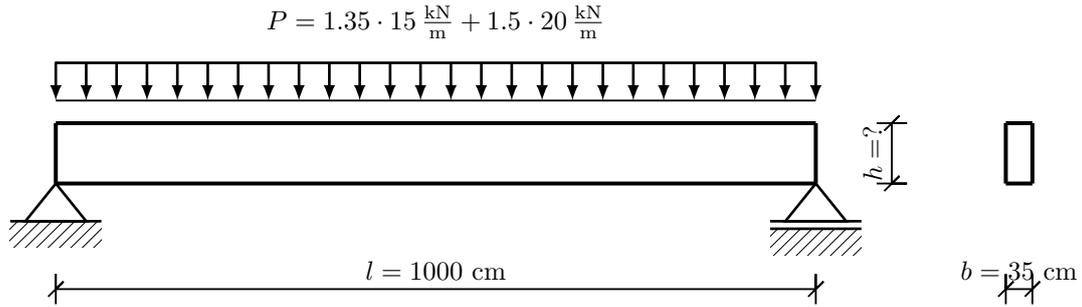
\begin{figure*}
\def\lx{10}
\def\ly{0.35}
\def\lz{0.8}
\centering
\begin{tikzpicture}

\point{LB}{0}{0};
\point{RB}{\lx}{0};
\point{RT}{\lx}{\lz};
\point{LT}{0}{\lz};

\point{MT}{0.5*\lx}{\lz+1}

\support{1}{LB};

\support{2}{RB};


\beam{2}{LB}{RB};
\beam{2}{RB}{RT};
\beam{2}{RT}{LT};
\beam{2}{LT}{LB};

\point{CLB}{\lx+2.5}{0};
\point{CRB}{\lx+2.5+\ly}{0};
\point{CRT}{\lx+2.5+\ly}{\lz};
\point{CLT}{\lx+2.5}{\lz};

\beam{2}{CLB}{CRB};
\beam{2}{CRB}{CRT};
\beam{2}{CRT}{CLT};
\beam{2}{CLT}{CLB};

\dimensioning{2}{RB}{RT}{+11.0}[$h=?$];
\dimensioning{1}{LB}{RB}{-1.4}[$l=1000\,\unit{cm}$];
\dimensioning{1}{CLB}{CRB}{-1.4}[$b=35\,\unit{cm}$];

\lineload{2}{LT}{RT}[0.5][0.5][0.04];

\notation{1}{MT}{$P=1.35\cdot15 \unit{\frac{kN}{m}} + 1.5\cdot20\unit{\frac{kN}{m}}$}[above];

\end{tikzpicture}
\caption{Geometry of the design problem of a bending beam with a constant distributed load (dead load and live load with safety factors of 1.35 and 1.5) and a rectangular cross section. The design variable, beam height is denoted by $h$}
\label{fig:design_problem}
\end{figure*}

In this work, a well-known example of a simply supported, reinforced, rectangular beam  has been chosen.
The design problem was originally published in \citep{everard1966reinforced} and illustrated in \reffig{fig:design_problem}.

It has been used to showcase different optimization schemes, e.g. \citep{Chakrabarty_1992}, \citep{Coello_1997}, \citep{Pierott_2021}. The objective is to reduce the overall GWP of the part.
This objective is particularly meaningful as the cement industry, accounts for approximately 8\% of the total anthropogenic GWP, \citep{Miller_2016}.
Reducing the environmental impact of concrete production becomes crucial in the pursuit of sustainable construction practices.
In addition, the reduction of the amount of cement in concrete is also correlated to the reduction of cost, as cement is generally the most expensive component of the concrete mix \citep{Paya_Zaforteza_2009}.
There are three direct ways to reduce the GWP of a given concrete part.
First, replace the cement with a substitute with a lower carbon footprint.
This usually changes mechanical properties and in particular, their temporal evolution.
Second, increase the amount of aggregates, therefore reducing the cement per volume.
This also changes effective properties and needs to be balanced with the workability and the limits due to the applications.
Third, decrease the overall volume of concrete, by improving the topology of the part.
In addition, when analyzing the whole life-cycle of a structure, both cost and GWP can be reduced by increasing the durability and therefore extending it's lifetime.
To showcase the proposed method's capability, two design variables have been chosen; the height of the beam and the ratio of ordinary Portland cement (OPC) to its replacement binder ground granulated blast furnace slag, a by-product of the iron industry.\\
In addition to the static design according to the standard, the problem is extended to include a key performance indicator related to the production process in a prefabrication factory that defines the time after which the removal of the formwork can be performed. To approximate this, the point in time when the beam can bear its own weight has been chosen a criterion. Reducing this time equates to being able to produce more parts with the same formwork.

\section{Methods}
\subsection{Design approches}\label{sec:design_approach}
\begin{figure}
\tikzstyle{ellipsenode} = [draw, ellipse, text width=1.3cm, minimum height=0.8cm, align=center, font=\footnotesize]
\tikzstyle{rectanglenode} = [draw, text width=1.5cm, minimum height=1cm, align=center, font=\footnotesize]
\tikzstyle{arrow} = [->]
\tikzstyle{arrowlabel} = [fill=white, fill opacity=0.75, text=black, inner sep=2pt, font=\footnotesize, midway]
\begin{tikzpicture}
  \node[rectanglenode] (load) at (0, 0) {loads};
  \node[rectanglenode] (geometry) at (0, -2) {geometry\\design};
  \node[ellipsenode] (structEng) at (2.5, 0) {structural\\engineer};
  \node[rectanglenode] (structDesign) at (2.5, -2) {structural\\design};
  \node[rectanglenode] (minMatProperty) at (2.5, -4) {minimial material properties};
  \node[ellipsenode] (matEng) at (5, 0) {material\\engineer};
  \node[rectanglenode] (matDesign) at (5, -2) {material design};
  
  \draw[arrow] (load) -- (structEng) node[arrowlabel] {input to};
  \draw[arrow] (geometry) -- (structEng) node[arrowlabel] {input to};
  \draw[arrow] (structEng) -- (structDesign) node[arrowlabel] {produces};
  \draw[arrow] (structDesign)--(minMatProperty)  node[arrowlabel] {includes};
  \draw[arrow] (minMatProperty)--(matEng)  node[arrowlabel] {input to};
  \draw[arrow] (matEng)--(matDesign)  node[arrowlabel] {produces};
\end{tikzpicture}
    \caption{Classical design approach, where the required minimal material properties are defined by the structural engineer which is then passed to the material engineer}
    \label{fig:standard_design}
\end{figure}
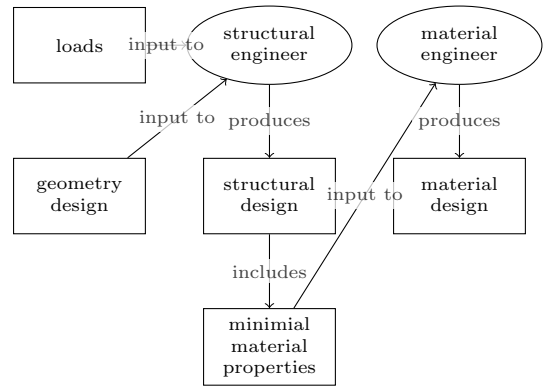

The conventional method of designing reinforced concrete structures is depicted in \reffig{fig:standard_design}. The structural engineer starts by chosing a suitable material (e.g. strength class C40/50) and designs the structure including the geometry (e.g. height of a beam) and the reinforcement. In the second step, this design is handed over to the material engineer with the constraint that the material properties assumed by the structural engineer have to be met. 
This lack of coordination strongly restricts the set of potential solutions since structural design and concrete mix design are strongly coupled, e.g. a lower strength can be compensated with a larger beam height.

\begin{figure}
\tikzstyle{ellipsenode} = [draw, ellipse, text width=1.3cm, minimum height=0.8cm, align=center, font=\footnotesize]
\tikzstyle{rectanglenode} = [draw, text width=1.5cm, minimum height=1cm, align=center, font=\footnotesize]
\tikzstyle{arrow} = [->]
\tikzstyle{arrowlabel} = [fill=white, fill opacity=0.75, text=black, inner sep=2pt, font=\footnotesize, midway]
\begin{tikzpicture}
  \node[rectanglenode] (matDesign) at (0, -2) {material design};
  \node[ellipsenode] (matEng) at (0, -4) {material\\engineer};
  \node[rectanglenode] (matProperty) at (0, -6) {material\\properties};
  \node[rectanglenode] (load) at (2.5, -2) {loads};
  \node[ellipsenode] (structEng) at (2.5, -4) {structural\\engineer};
  \node[rectanglenode] (structProp) at (2.5, -6) {structural\\design};
  \node[rectanglenode] (geometry) at (5, -2) {geometry\\design};
  
  \draw[arrow] (matDesign) -- (matEng) node[arrowlabel] {input to};
  \draw[arrow] (load) -- (structEng) node[arrowlabel] {input to};
  \draw[arrow] (geometry) -- (structEng) node[arrowlabel] {input to};
  \draw[arrow] (matEng) -- (matProperty) node[arrowlabel] {predicts};
  \draw[arrow] (structEng) -- (structProp) node[arrowlabel] {predicts};
  \draw[arrow] (matProperty) -- (structEng) node[arrowlabel] {input to};
\end{tikzpicture}
    \caption{Proposed inverse-design approach that is integrated into a holistic optimization approach}
    \label{fig:proposed_workflow}
\end{figure}
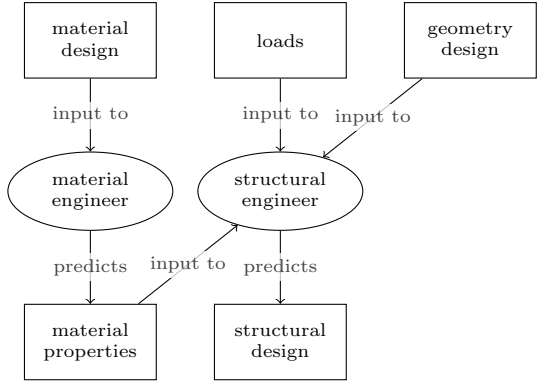

An alternative design workflow is illustrated in \reffig{fig:proposed_workflow} which entails inverting the classical design pipeline.
The material composition is the input to the material engineer who predicts the corresponding mechanical properties of the material. This includes Key Performance Indicators (KPIs) related to the material, e.g. viscosity/slump test, or simply the water/cement ratio. In a second step, a  structural analysis is performed with the material properties as input. This step outputs the structural KPIs such as the load bearing capacity, the expected lifetime (for a structure susceptible to fatigue) or the GWP of the complete structure. These two (coupled) modules are used within an optimization procedure to estimate the optimal set of input parameters (both on the material level as well as on the structural level). One of the KPIs is chosen as the objective function (e.g. GWP) and others as constraints (e.g. load-bearing capacity larger than the load, cement content larger than a threshold, viscosity according to the slump test within a certain interval). Note that in order to use such an inverse-design procedure, the forward modeling workflow needs to be automated and subsequently the information needs to be  efficiently back-propagated. 

\begin{figure*}
\tikzstyle{ellipsenode} = [draw, text width=2.0cm, minimum height=1.5cm, align=center, font=\footnotesize, rounded corners=0.5cm, thick, fill=gray!50!white]
\tikzstyle{rectanglenode} = [draw, text width=1.6cm, minimum height=1cm, align=center, font=\footnotesize, fill=gray!30!white]
\tikzstyle{kpi} = [draw, ellipse, text width=1.8cm, minimum height=1.5cm, align=center, font=\footnotesize, fill=gray!7!white]
\tikzstyle{arrow} = [->]
\begin{tikzpicture}
  \draw [draw=none, fill=gray!10!white] (-5.7,4.5) rectangle (-1.1,1.4);
  \node[rectanglenode] (parameter) at (-2.2, 3.6) {input/output};  
  \node[ellipsenode] (function) at (-4.4, 3.6) {function/model};
  \node[kpi] (parameter) at (-3.3, 2.1) {KPI};
  \node[rectanglenode] (geometry) at (-2.2, 0.0) {geometry};  
  \node[rectanglenode] (loads) at (0.0, 1.9) {loads};
  \node[rectanglenode] (mix) at (3.3, 3.6) {mix composition};
  \node[ellipsenode] (ML model mech) at (4.7, 1.8) {ML model paste mech. (Sec. \ref{sec:model_learning})};
  \node[ellipsenode] (ML model hyd) at (7.7, 1.8) {ML model paste hyd. (Sec. \ref{sec:model_learning})};
  \node[rectanglenode] (steel properties) at (-4.4, 0) {steel\\properties};
  \node[rectanglenode] (aggregate properties) at (2.2, 0) {aggregate\\properties};
  \node[rectanglenode] (pasteOther) at (4.4, 0) {paste prop. (excl. $E_{28}$, $f_{c,28}$)};
  \node[rectanglenode] (pasteE28) at (6.6, 0.0) {paste\\$E_{28}$, $f_{c,28}$};
  \node[rectanglenode] (concrete) at (4.4, -3.6) {concrete properties};
  \node[rectanglenode] (hydration parameters) at (7.7, -3.6) {parameters hydration model};
  \node[ellipsenode] (GWPcomputation) at (-4.4, -3.6) {computation Global Warming Potential (Sec. \ref{sec:GWP})};
  \node[ellipsenode] (homogenization) at (3.3, -1.8) {homogenization (Sec. \ref{sec:homogenization_model})};
  \node[ellipsenode] (designEC) at (0.0, -5.4) {design\\(Eurocode) (Sec \ref{sec:beam_design})};
  \node[ellipsenode] (FEM) at (3.3, -5.4) {FEM solver (Sec \ref{sec:FE_model})};
  \node[rectanglenode] (reinforcement) at (-2.7, -5.4) {amount of  reinforcement};
  \node[kpi] (GWP) at (-4.4, -7.4) {Global Warming Potential};
  \node[kpi] (maxReinforcement) at (0.0, -7.4) {reinforcement geometry constraint};
  \node[kpi] (formwork) at (3.3, -7.4) {time of demoulding};
  \node[kpi] (maxT) at (6.6, -7.4) {max. temperature during hydration};
  \draw[arrow] (loads) -- (designEC) ;
  \draw[arrow] (geometry) -- (designEC) ;
  \draw[arrow] (steel properties) -- (designEC);
  \draw[arrow] (mix) -- (homogenization) ;
  \draw[arrow] (pasteOther) -- (homogenization);
  \draw[arrow] (pasteE28) -- (homogenization);
  \draw[arrow] (aggregate properties) -- (homogenization);
  \draw[arrow] (homogenization) -- (concrete) ;
  \draw[arrow] (concrete) -- (designEC);
  \draw[arrow] (concrete) -- (FEM);
  \draw[arrow] (designEC) -- (reinforcement) ;
  \draw[arrow] (reinforcement) -- (GWPcomputation);
  \draw[arrow] (designEC) -- (maxReinforcement) ;
  \draw[arrow] (geometry) -- (GWPcomputation) ;
  \draw[arrow] (aggregate properties) -- (GWPcomputation) ;
  \draw[arrow] (steel properties) -- (GWPcomputation);
  \draw[arrow] (mix) -- (GWPcomputation) ;
  \draw[arrow] (mix) -- (ML model hyd) ;
  \draw[arrow] (mix) -- (ML model mech) ;
  \draw[arrow] (ML model hyd) -- (hydration parameters) ;  
  \draw[arrow] (ML model mech) -- (pasteE28) ;
  \draw[arrow] (hydration parameters) -- (FEM) ;
  \draw[arrow] (GWPcomputation) -- (GWP) ;
  \draw[arrow] (loads) -- (FEM) ;
  \draw[arrow] (steel properties) -- (FEM) ;
  \draw[arrow] (geometry) -- (FEM) ;
  \draw[arrow] (FEM) -- (formwork) ;
  \draw[arrow] (FEM) -- (maxT) ;
  \draw[arrow] (pasteOther) -- (GWPcomputation) ;
\end{tikzpicture}
    \caption{Workflow to compute KPIs from input parameters}\label{fig:snakemake_workflow}
\end{figure*}
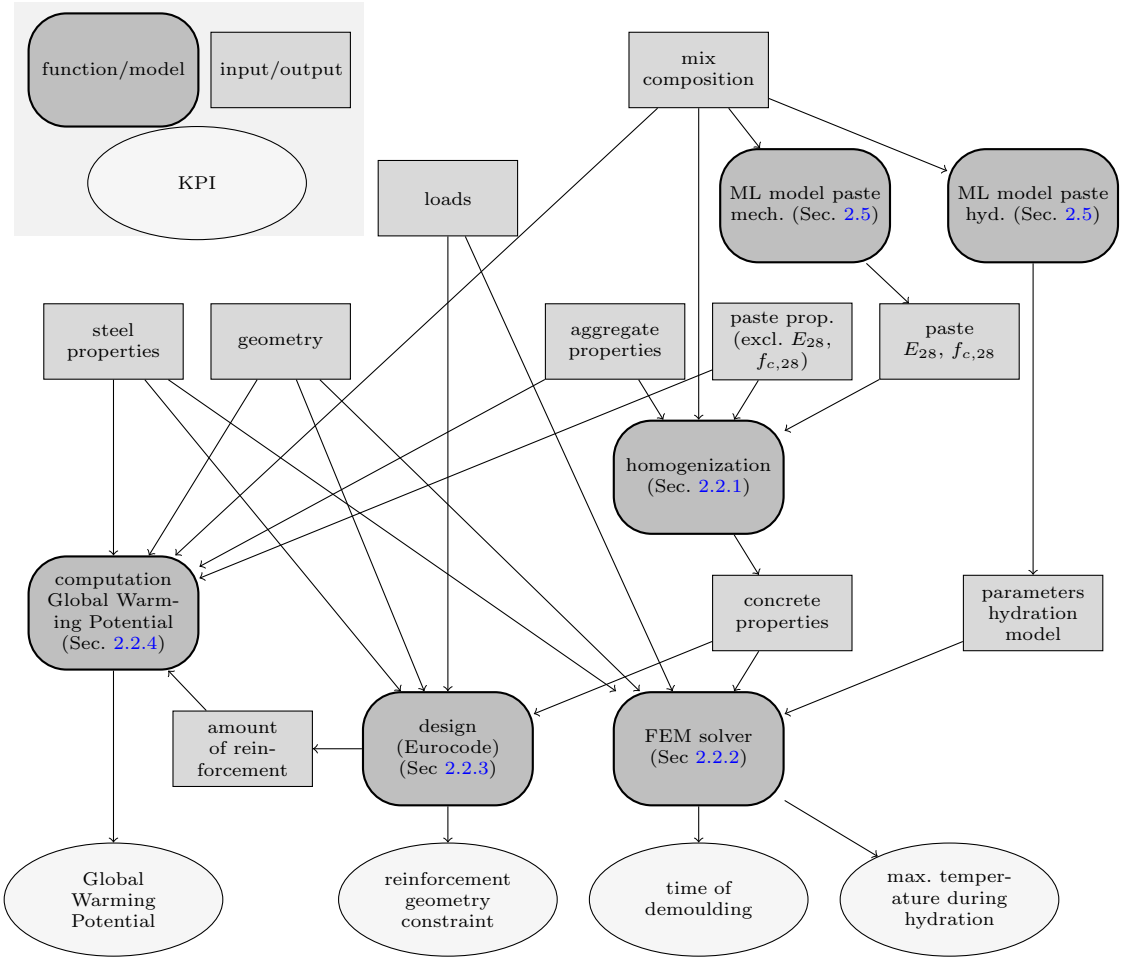

The paper aims to present the proposed methodological framework as well as illustrate its capabilities in the design of a precast concrete element with the objective of  reducing  the GWP. The  constraints employed are  related to the structural performance after 28 days as well as the maximum time of demoulding after 10 hours. The design/optimization variables are,  on the structural level, the height of the beam, and on the material level the composition of the binder as a mixture of Portland cement and slag. 
The complete workflow is illustrated in \reffig{fig:snakemake_workflow}. 
\subsection{Workflow for predicting key performance indicators}\label{sec:models}
The workflow consists of four major steps. In a first step, the cement composition (blended cement and slag) defined in the mix composition is used to predict the mechanical properties of the cement paste. This is done using a data-driven approach  as  discussed in \refsec{sec:calibration}.
In a second step, homogenization is used in order to compute the  effective, concrete properties based on cement paste and aggregate data.
An analytical function is applied for the homogenization based on the Mori-Tanaka scheme \cite{mor_1973_asi}. 
The third step involves  a multi-physics, finite element model with two complex constitutive models - a hydration model, which computes the evolution of the degree of hydration, considering the local temperature and the heat released during the reaction and a  mechanical model which simulates the temporal evolution of the mechanical properties   assuming that those 
 depend on the degree of hydration. 
The fourth and last model is based on a design code to estimate the amount of reinforcement and predict the load bearing capacity after 28 days.
Subsequent sections will provide insights into how these models function within the optimization framework.

\subsubsection{Homogenized Concrete Parameters}\label{sec:homogenization_model}
Experimental data for estimating the compressive strength is obtained from concrete specimens measuring the homogenized response of cement paste and aggregates. The mechanical properties of aggregates are known, whereas the cement paste properties have to be inversely estimated. The calorimetry is directly performed for cement paste.

In order to relate macroscopic mechanical properties to the individual constituents (cement paste and aggregates), an analytical homogenization procedure is used. The homogenized effective concrete properties are the Young's modulus $\eMod$, the Poisson's ratio $\poission$, the compressive strength $\fc$, the density $\density$, the thermal conductivity $\thermCond$, the heat capacity $\heatCapSpecific$ and the total heat release $\heatInf$.
Depending on the physical meaning, these properties need slightly different methods to estimate the effective concrete properties.
The elastic, isotropic properties $\eMod$ and $\poission$ of the concrete are approximated using the Mori-Tanaka homogenization scheme \citep{mor_1973_asi}.
The method assumes spherical inclusions in an infinite matrix and considers the interactions of multiple inclusions. Details given in \ref{ssec:mt_elastic}.\\
The estimation of the concrete compressive strength $\fcEff$ follows the ideas of \citep{nev_2018_mcam}. The premise is that a failure in the cement paste will cause the concrete to crack. The approach is based on two main assumptions.
First, the Mori-Tanaka method is used to estimate the average stress within the matrix material $\stressMatrix$. Second, the von Mises failure criterion of the average matrix stress is used to estimate the uniaxial compressive strength (see \ref{ssec:compressivestrength}).

\begin{figure*}[!htbp]%
	\centering
	\includegraphics[width=\textwidth]{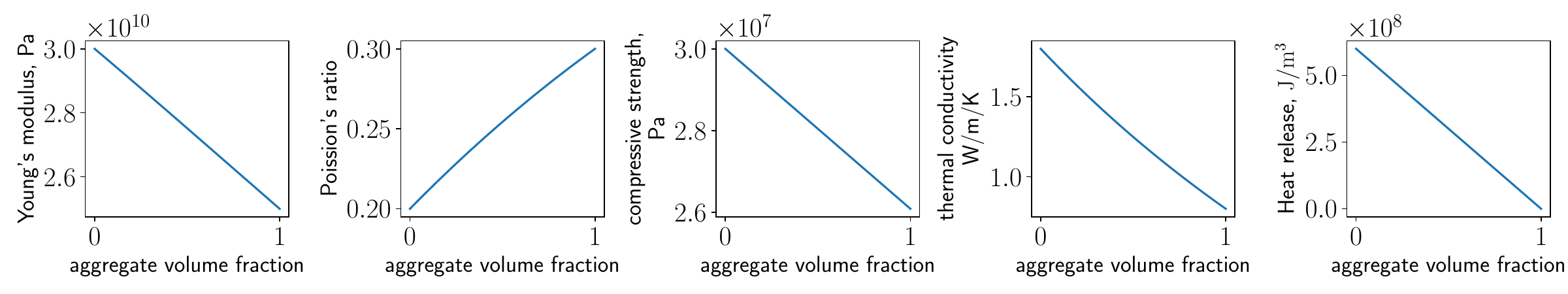}
	\caption{\label{fig:homogenization}Influence of aggregate ratio on effective concrete properties}
\end{figure*}

\mbox{Table \ref{tab:homogenizationproperties}} gives an overview of the material properties of the constituents used in the subsequent sensitivity studies.
\begin{table}[b]
	\begin{center}
		\begin{minipage}{\columnwidth}
			\caption{Properties of the cement paste and aggregates used in subsequent sensitivity studies}\label{tab:homogenizationproperties}
			\setlength{\tabcolsep}{3pt}
			\begin{tabular}{lccccccc}
				\toprule
				Phase & $\eMod$ & $\poission$ & $\fc$ & $\density$ & $\thermCond$ & $\heatCapSpecific$ & $\heatInf$\\
				  & $\left[\unit{\pasteEunit}\right]$  & $\left[-\right]$  & $\left[\unit{\pastefcunit}\right]$  & $\left[\unit{\pasterhounit}\right]$  & $\left[\unit{\pasteCunit}\right]$  & $\left[\unit{\pastekappaunit}\right]$  &  $\left[\unit{\pasteQunit}\right]$ \\
				\midrule
			Paste	& \pasteE & \pastenu & \pastefc & \pasterho & \pasteC & \pastekappa &  \pasteQ \\
			Aggr.	& \aggregatesE & \aggregatesnu & - & \aggregatesrho & \aggregatesC & \aggregateskappa &  0 \\
				\botrule
			\end{tabular}
		\end{minipage}
	\end{center}
		
\end{table}
The effective properties as a function of the aggregate content are plotted in \reffig{fig:homogenization}. Note that both extremes (0 - pure cement and 1 - only aggregates) are purely theoretical.

For the considered example, the relations are close to linear. This can change, when the difference between the matrix and the inclusion properties is more pronounced or more complex micro mechanical mechanisms are incorporated, as air pores or the interfacial transition zone. Though not done in this paper, these could be considered within the chosen homogenization scheme by adding additional phases, c.f. \citep{nee_2012_ammf}.
Homogenization of the thermal conductivity is also based on the Mori-Tanaka method, following the ideas of \citep{str_2011_mbeo} with details given in appendix \ref{ssec:thermalconductivity}.
The density $\density$, the heat capacity $\heatCapSpecific$ and the total heat release $\heatInf$ can be directly computed based on their volume average.
As example for the volume averaged quantities, the heat release is shown in \reffig{fig:homogenization} as it exemplifies the expected linear relation of the volume average as well as the zero heat output of a theoretical pure aggregate.

\subsubsection{Hydration and evolution of mechanical properties}\label{sec:FE_model}
Due to a chemical reaction (hydration) of the binder with water, Calcium-silicate hydrates (CSH) form that lead to a temporal evolution of concrete strength and stiffness. The reaction is exothermal and the kinetics are sensitive to the temperature. The primary model simulates the hydration process and computes the temperature field $\temp$ and the degree of hydration (DOH) $\DOH$ (see Eq. (\ref{eq:heat1}, \ref{eq:doh}) in the appendix). The latter characterizes the degree of hydration that condenses the complex chemical reactions into a single scalar variable.
The thermal model depends on three material properties, the effective thermal conductivity $\thermCondEff$, the specific heat capacity $C$ and the heat release $\dQdt$.
The latter is governed by the hydration model, characterized by six parameters:
$\hydParBone, \hydParBtwo, \hydParEta, \tempRef, \activE$ and $\DOHmax$.
The first three $\hydParBone, \hydParBtwo$ and $\hydParEta$ are parameters characterizing the shape of the evolution of the heat release. $\tempRef$ is the reference temperature for which the first three parameters are calibrated (Based on the difference between the actual and the reference temperature, the heat released is scaled). The sensitivity to the temperature is characterized by the activation energy $\activE$. $\DOHmax$ is the maximum degree of hydration that can be reached. Following \citep{Mills1966fico}, the maximum degree of hydration is estimated based on the water to binder ratio $\wc$, as $\DOHmax = \frac{1.031\,\wc}{0.194 + \wc}$.

By assuming the DOH to be the fraction of the currently released heat with respect to its theoretical potential $\heatInf$, the current degree of hydration is estimated as $\DOH(\zeit) = \frac{\heat(\zeit)}{\heatInf}$.
As the potential heat release is also difficult to measure as it takes a long time to fully hydrate and will only do so under perfect conditions, we identify it as an additional parameter in the model parameter estimation.
For a detailed model description see Appendix \ref{appendix:fem}.
In addition to influencing the reaction speed, the computed temperature is used to verify that the maximum temperature during hydration does not exceed a limit of $\tempLimit = \inputtemperaturelimit$\textdegree C. 
Above a certain temperature, the hydration reaction changes (e.g. secondary ettringite formation) and, additionally, volumetric changes in the cooling phase correlate with cracking and reduced mechanical properties.
The maximum temperature is implemented as a constraint for the optimization problem (see \refeq{eq:concstraintT}).\\
The evolution of the Young's modulus $\eMod$ of a linear-elastic material model is modelled  as a function of the degree of hydration (details in \refeq{eq:EwrtDOH}). In a similar way, the compressive strength evolution is computed (see \refeq{eq:fcwrtDOH}), which is utilized to determine a failure criterion based on the computed local stresses \refeq{eq:constraintStress} related to the time when the formwork can be removed.
For a detailed description of the parameter evolution as a function of the degree of hydration see Appendix \ref{appendix:fem_evolution}.
\reffig{fig:parameterEvolution} shows the influence of the different parameters.
In addition to the formulations given in \citep{car_2016_mamt} which depend on a theoretical value of parameters for fully hydrated concrete at $\DOH = 1$, this work reformulates the equations, to depend on the 28 day values $\eModTwentyEight$ and $\fcTwentyEight$ as well as the corresponding $\DOHTwentyEight$ which is obtained via a simulation.
This allows to directly use the experimental values as input.
In \reffig{fig:parameterEvolution}, $\DOHTwentyEight$ is set to $\evoExAlphaTx$.\\
\begin{figure*}[!htpb]%
	\centering
	\includegraphics[width=0.8\textwidth]{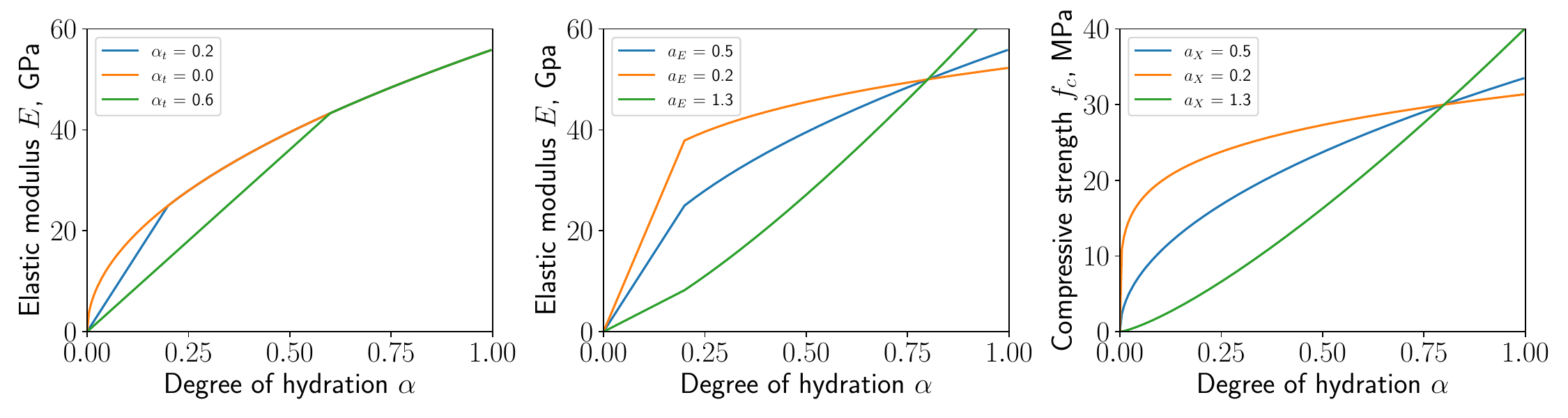}
	\caption{Influence of parameters $\DOHt, \stiffExp$ and $\strengthCExp$ on the evolution the Young's modulus and the compressive strength with respect to the degree of hydration $\DOH$.
	Parameters: $\eModTwentyEight = \evoExE \unit{GPa}$, $\stiffExp = \evoExaE$, $\DOHt = \evoExAlphaT$, $\strengthCExp = \evoExafc$,  $\fcTwentyEight = \evoExfc$ $\unit{N/mm^2}$, $\DOHTwentyEight = \evoExAlphaTx$.}
 \label{fig:parameterEvolution} 
\end{figure*}

\subsubsection{Beam design according to EC2}\label{sec:beam_design}
\begin{figure*}[ht]%
	\centering
	\includegraphics[width=0.8\textwidth]{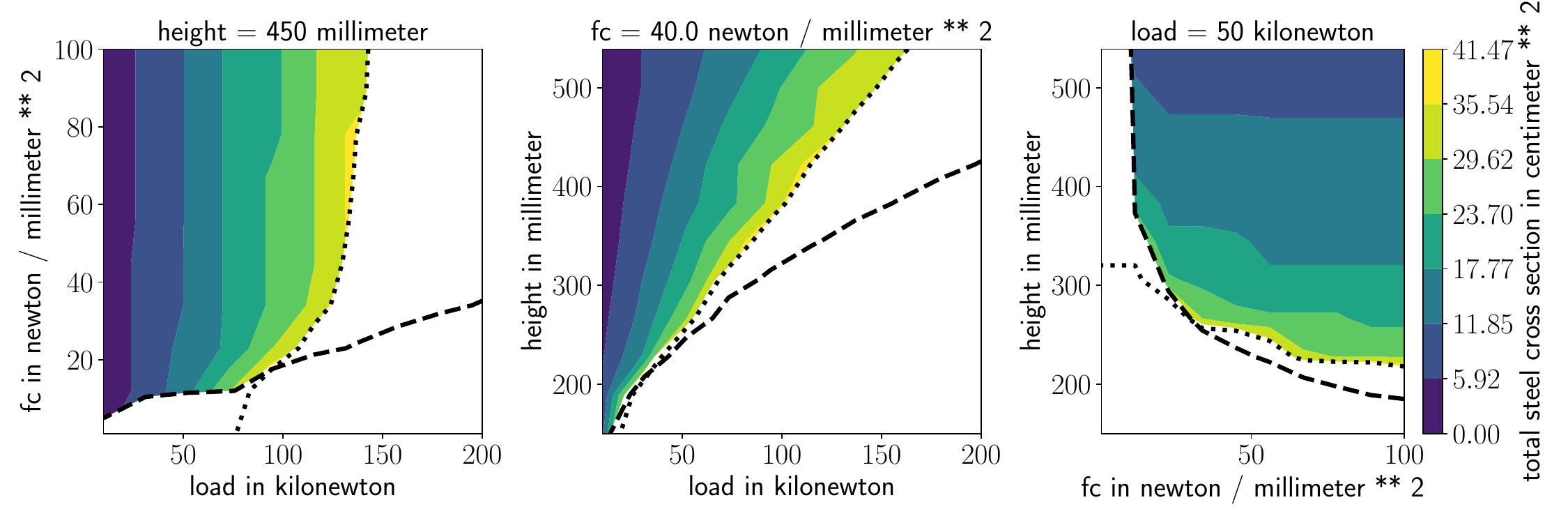}
	\caption{Influence of beam height, concrete compressive strength and load in the center of the beam on the required steel. The dashed lines represent the minimum compressive strength constraint \refeq{eq:constraintfc}, the dotted lines the geometrical constraint from the spacing of the bars \refeq{eq:constraintGeo} \label{fig:beamdesign}}
\end{figure*}

The design of the reinforcement and the computation of the load-bearing capacity is performed based on \citep{DIN1992-1-1} according to \refeq{eq:Areq} with a detailed explanation in the appendix. To ensure that the design is realistic, the continuous cross section is transformed into a discrete number of bars with a diameter chosen from a list.
This is visible in \reffig{fig:beamdesign} by the step-wise increase in cross sections.
The admissible results are restricted by two constrains.
One is coming from a minimal required compressive strength \refeq{eq:constraintfc}, visualized as dashed line.
The other, based on the available space to place bars with admissible spacing \refeq{eq:constraintGeo}, visualized as the dotted line.
Further detail on the computation are given in Appendix \ref{appendix:beam}.
A sensitivity study for the mutual interaction and the constraints is visualized in \reffig{fig:beamdesign}. The parameters for the sensitivity study are given in Table \ref{tab:beamdesigninput}.

\subsubsection{Computation of GWP}\label{sec:GWP}
The computation of the global warming potential is performed by multiplying the volume content of each individual material by its specific global warming potential. The values used in this study are extracted from \citep{Braga2017} and listed in Table \ref{table:GWP}.
\begin{table}[bt]
    \centering
    \begin{tabular}{cc}
         \hline
         material& GWP $\left[\unit{\dfrac{kg_{CO_2 eq}}{m^3}}\right]$  \\
         \hline
         portland cement& 0.95 \\
         slag& 0.18 \\
         aggregates& 0.025 \\
         water& 0.000133 \\
         steel& 1.42 \\
         \hline
    \end{tabular}
    \caption{Specific global warming potential of the raw materials used in the design.}
    \label{table:GWP}
\end{table}

The values are certaintly a great source of discussion in the community and serve here only a exemplary values. This is due to the question, what exactly to include in the GWP computation, e.g. the transport of materials is difficult to generally include, there are always local conditions (e.g. the GWP of the energy sources used in the cement production depends on the amount of green energy in that country), the time span (complete life cycle analysis vs. production) is a point of debate and finally the usage of by-products (slag is currently a by-product of steel manufacturing and thus its GWP is considered to be small).

\subsection{Experimental Data}\label{sec:exp_data}
This section describes the data used for learning the missing (probabilistic) links (detailed in \refsec{sec:model_learning}) between the slag-binder mass ratio $r_{sb}$ and physical model parameters. The slag-binder mass ratio $r_{sb}$ is the mass ratio between the amount of Blast Furnace Slag and the binder (sum of  Blast Furnace Slag (BFS) and Ordinary Portland Cement (OPC)). The data is sourced from \citep{gruyaert2011}. In particular, we are concerned about the parameter estimation for the concrete homogenization discussed in \refsec{sec:homogenization_model} and the hydration model in \refsec{sec:FE_model}. 

For concrete homogenization, six different tests for varying ratios $r_{sb} = \{0.0,0.15,0.3,0.5,0.7,0.85\}$ are available for the concrete compressive strength $f_{c}$ after 28 days. For the concrete hydration, we utilize isothermal calorimetry data at $20^{\circ}C$. We have temporal evolution data of cumulative heat of hydration $\hat{\bm Q}$ for four different values of $r_{sb} = \{0.0,0.30,0.50,0.85\}$, as illustrated in \reffig{fig:prediction_hyd_model}.
For details on other material parameters and phenomenological values used to obtain the data, the reader is directed to \citep{gruyaert2011}.

\subsubsection{Young's modulus \texorpdfstring{$\eMod$}{E} based on \texorpdfstring{$fc$}{fc}}
The dataset does not encompass information about the Young's modulus. 
Given its significance for  the FEM simulation, we resort to a phenomenological approximation derived from \citep{ACI363}. 
This approximation relies on the compressive strength $\fc$ and the density $\density$ to estimate the Young's modulus
\begin{align}
	\eMod = 3320 \sqrt{\fc} + 6895 \left( \dfrac{\density}{2320}\right)^{1.5},
\end{align}
with $\density$ in $\unit{\dfrac{kg}{m^3}}$, $\fc$ and $\eMod$ in $\unit{MPa}$.

\subsection{Model learning and optimization}
\label{sec:calibration}

\begin{figure}[!htpb]
	\centering
	\begin{tikzpicture}
		\node (theta_1) at (-1,0) {$\theta_1$};
		\node (theta_2) at (-1,-2) {$\theta_2$};
		\node (x_1) at (0,0) [circle, draw, dotted] {$x_1$};
		\node (x_2) at (0,-2) [circle, draw, dotted] {$x_2$};
		\node (b) at (1,0) [circle, draw] {$\bm{b}$};
		\node (y_1) at (3,1) [rectangle, draw] {$y_{c_1}$};
		\node (y_2) at (3,0) [rectangle, draw] {$y_{c_2}$}; 
		\node (y_3) at (3,-1) [rectangle, draw] {$y_{c_3}$}; 
		\node (y_4) at (3,-2) [rectangle, draw] {$y_o$}; 
		\node (C_1) at (5,1) [rectangle, draw] {$\mathcal{C}_1\left(y_{c_1}(\cdot)\right)$}; 
		\node (C_2) at (5,0) [rectangle, draw] {$\mathcal{C}_2\left(y_{c_2}(\cdot)\right)$}; 
		\node (C_3) at (5,-1) [rectangle, draw] {$\mathcal{C}_3\left(y_{c_3}(\cdot)\right)$}; 
		\node (O) at (5,-2) [rectangle, draw] {$\mathcal{O}\left(y_o(\cdot)\right)$}; 
		
		\graph {
			
			(x_1) -> {(b),(y_4)};
			(b) -> {(y_1),(y_2),(y_3)};
			(y_1) -> (C_1);
			(y_2) -> (C_2);
			(y_3) -> (C_3);
			(y_4) -> (O);
			(theta_2) -> (x_2);
			(theta_1) -> (x_1);
			(x_2) -> {(y_1),(y_2),(y_3),(y_4)};
			
		};
	\end{tikzpicture}
	\caption{\emph{Stochastic computational graph for the constrained optimization problem for  performance-based concrete design:} The circles represent \textit{stochastic nodes} and  rectangles  \textit{deterministic nodes}.
	The design variables are denoted by $\bx = (x_1,x_2)$. The vector $\bm{b}$ represents the  unobserved model parameters  which are needed in order to link the  KPIs $\bm{y} = (y_o,\{y_{c_i}\}_{i=1}^I)$ with the design variables $\bx$. Here  $y_o$ represents the model-output appearing in the optimization objective and $y_{c_i}$ represents the model output appearing in  the $i^{th}$ constraint. The objective function is given by $\calO$ and the $i^{th}$ constraint by $\calC_i$. They are not differentiable w.r.t to $x_1,x_2$ (Hence $x_1$ and $x_2$ is dotted). The variables $\bt$ are auxiliary and are used in the context of Variational Optimization   discussed in \refsec{sec:non-diff_obj_cons}. 
	Several other deterministic nodes are present between the random variables $\bm{b}$ and the KPIs $\bm{y}$ but they are omitted for brevity. The physical meaning of the variables used is detailed in Table \ref{table:stoc_graph_variables}}
	\label{fig:stochastic graph demonstrator}
\end{figure}
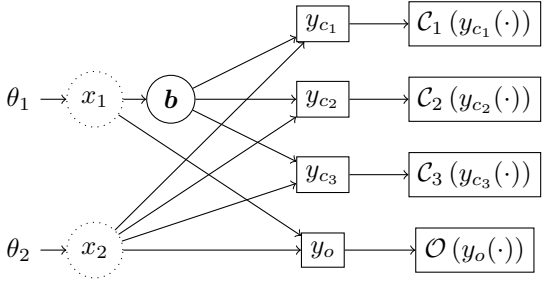

The workflow illustrated in \reffig{fig:snakemake_workflow}, which builds the link between the parameters relevant to the concrete mix design and the KPIs involving the environmental impact and the structural performance can be represented in terms of the probabilistic graph shown  in \reffig{fig:stochastic graph demonstrator}. As  discussed in the  Introduction (section  \ref{sec:introduction}), the goal of the present study is to find the value of the design variables $\bx$ (concrete mix design, beam geometry) which minimizes the objective $\calO$ (environmental impact), while satisfying a given set of constraints $\calC_i$ (beam design criterion, structural performance etc.). This necessitates forward and backward information flow in the presented graph. 
The forward information flow is necessary to compute the KPIs for given values of the design variables and the backward information is essentially the sensitivities of the objective and the constraints with respect to the design variables that  enable gradient-based optimization. Establishing the information flow poses challenges, which we attempt to tackle with the methods proposed as follows:
\begin{itemize}
	\item \emph{Data-based model learning:} The physics-based  models discussed in  (\refsec{sec:homogenization_model} and  \refsec{sec:FE_model}) are used to compute various KPIs (discussed in \reffig{fig:stochastic graph demonstrator}). These depend on some model parameters denoted by $\bs{b}$  which are unobserved (latent) in the experiments performed. The model parameters need not only be inferred on the basis of experimental data but also their dependence on the design variables $\bx$ is required in order to be integrated in the optimization framework. In addition, the noise in the data (aleatoric) or the incompleteness of data (epistemic) introduce uncertainty.  
To this end we propose learning {\em probabilistic} links by employing experimental data as  discussed in detail in \refsec{sec:model_learning}. 
	\item \emph{Optimization under uncertainty:}
	The aforementioned uncertainties as well as additional randomness that might be present in the associated links necessitate reformulating the optimization problem (i.e. objectives/constraints) as one of optimization under uncertainty. In turn this gives rise to new challenges in order to compute the needed derivatives of the KPIs with respect to the design variables  which are discussed in  \refsec{sec:optimization}. 
\end{itemize}

\begin{table}[!htbp]
	\centering
	\begin{tabular}{l p{6cm}}
            \hline
            Variables & Physical meaning \\
            
			\hline
			$x_1$ & Mass ratio of Blast Furnace Slag (BFS) and Ordinanry portland cement (OPC) $r_{sb}$\\
			$x_2$ & Beam height $h$ \\
			$\bm b$ & Vector of the input, model parameters to the homogenization and hydration model (\refsec{sec:homogenization_model} and \refsec{sec:FE_model} respectively)\\
			$y_{c_1}$ &  Required steel reinforcement area (Eq. \ref{eq:Areq}) \\
			$\mathcal{C}_1\left(y_{c_1}(\cdot)\right)$ & Beam design constraint (Eq. \ref{eq:design_cons} )\\
			$y_{c_2}$ &  Max. temperature reached (Sec. \ref{sec:appendix_constraints})  \\
			$\mathcal{C}_2\left(y_{c_2}(\cdot)\right)$ & Temperature constraint (Eq. \ref{eq:concstraintT})\\
			$y_{c_3}$ &  Time of demolding (Sec. \ref{sec:appendix_constraints}) \\
			$\mathcal{C}_3\left(y_{c_3}(\cdot)\right)$ & Time constraint based on yield strength (Eq. \ref{eq:constraintStress})\\
			$y_{o}$ &  The GWP of the beam (Sec. \ref{sec:GWP}) \\
			$\mathcal{O}_o\left(y_{o}(\cdot)\right)$ & Objective corresponding to the beam GWP.\\
			\hline
		\end{tabular}
	\caption{Physical meaning of the variables used in \reffig{fig:stochastic graph demonstrator}}
	\label{table:stoc_graph_variables}
\end{table}

\subsection{Probabilistic links based on data and physical models}\label{sec:model_learning}
\begin{figure}[!htpb]
	\centering
	\includegraphics[width=0.35\textwidth]{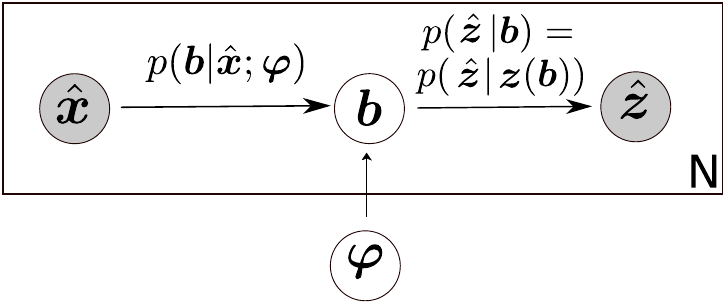}
	\caption{\emph{Probabilistic graph for the data and physical model based model learning:} the shaded nodes are the observed and unshaded are the unobserved (latent) nodes}
	\label{fig:cal_graph_general}
\end{figure}
This section deals with learning a (probabilistic) model linking the design variables $\bx$ and the input parameters of the physics-based  models i.e. concrete hydration and concrete homogenization. 
A graphical representation is contained in \reffig{fig:cal_graph_general}.
Therein,  $\{\hat{\bx}^{(i)},  \hat{\bz}^{(i)}\}_{i=1}^N$ denote the observed data-pairs and $\bm b$ denotes a vector of unknown and unobserved  parameters of the physics-based models and $\bz(\bm b)$ the model outputs. The latter relate to an experimental observation $\hat{\bz}^{(i)}$ as  $\hat{\bz}^{(i)}=\bz(\bm b^{(i)})+noise $ which gives rise to a likelihood $p(\hat{\bz}^{(i)} \mid \bm{z}(\bm b^{(i)}))$. 
We further postulate a probabilistic relation between  $\hat{\bx}$ and $\bm b$ that expressed by the conditional $p(\bs{b}\mid \bx; \varphi)$  which depends on unknown parameters $\bm \varphi$. 
The physical meaning of the aforementionned variables and model links as  well as of the relevant data is presented  in Table \ref{table:model_learning}. The elements introduced  above suggest a Bayesian formulation which can quantify inferential uncertainties in the unknown parameters and propagate it in the model predictions  \citep{koutsourelakis2016big}, as detailed in the next section. 
\begin{table}[!htbp]
	\centering
	\begin{tabular}{p{2.5cm} p{2.5cm} p{1.5cm}} 
			\hline
			$\bm b$ (model input) & $\hat{\bz}$ (observed data) & $\bz(\bm b)$\newline(physics-based model) \\
			\hline 
			hydration model input parameters (\refsec{sec:FE_model}) $\bb =[B_1, B_2, \eta, Q_{pot}]$ & Heat of hydration $\bm Q$ & concrete hydration model\\
            \\
			cement paste compressive strength ($f_{c,paste}$), cement paste Youngs Modulus ($E_{paste}$) (\refsec{sec:homogenization_model}) & concrete compressive strength ($f_{c}$), concrete Youngs Modulus ($E_{c}$)  & concrete homogenization model\\
			\hline
		\end{tabular}
	\caption{Physical meaning of the variables/links used in \reffig{fig:cal_graph_general}. Note $\hat{\bx}$ is the slag-binder mass ratio $r_{sb}$ for the all the cases presented above. }
	\label{table:model_learning}
\end{table}

\subsubsection{Expectation-Maximization}\label{sec:EM}

Given  $N$ data-pairs $\mathcal{D}_N=\{ \hat{\bx}^{(i)},  \hat{\bz}^{(i)}\}_{i=1}^N$ consisting of different concrete mixes and corresponding outputs, we would like to  to infer the corresponding $\bs{b}^{(i)}$, but more importantly the relation between $\bx$ and $\bs{b}$  which would be of relevance for downstream, optimization tasks discussed in \refsec{sec:optimization}.

We postulate a probabilistic relation between $\bx$ and $\bs{b}$ in the form of a conditional density $p(\bs{b}\mid \bx; ~\bm{\varphi})$ parametrized by  $\bm \varphi$. E.g.:
\begin{align}\label{eq:prior}
	p(\bs{b}\mid \bx; ~\bm{\varphi})=\mathcal{N}(\bs{b}\mid ~\bs{f}_{\varphi} (\bx), \bs{S}_{\varphi}(\bx))
\end{align}
where  $\bs{f}_{\varphi} (\bx)$ represents a fully connected, feed-forward neural network parametrized by $\bm{w}$ (further details discussed in \refsec{sec:numericalexperiments}), $\bs{S}_{\varphi}(\bx)=\bs{LL}^T$ denotes the covariance matrix where $\bs{L}$ is lower-triangular. Hence the parameters $\bm \varphi$ to be learned correspond to 
$\bm{\varphi} =\{\bm{w}, \bm{L}\}$. 
We assume that the observations $\hat{\bz}^{(i)}$ are contaminated with Gaussian noise, which gives rise to the likelihood:
\begin{align}\label{eq:likelihood}
    p(\hat{\bz}^{(i)} \mid \bz(\bs{b}^{(i)})) = \mathcal{N} (\hat{\bz}^{(i)} \mid \bz(\bs{b}^{(i)}), \bm{\Sigma}_{\ell}).
\end{align}
The covariance  $\bm{\Sigma}_{\ell}$  depends on the data used and is discussed in \refsec{sec:numericalexperiments}.

Given \refeq{eq:prior} and \refeq{eq:likelihood}, one can observe that $\bs{b}^{(i)}$ (i.e. the unobserved model inputs for each concrete mix $i$) and  $\bm \varphi$ would need to be inferred simultaneously. 
In the following we obtain point-estimates  $\bm \varphi^*$  for  $\bm \varphi$, by maximizing the marginal log-likelihood $p(\mathcal{D}_N \mid \bm \varphi)$  (also known as log-evidence) i.e. the probability that the observed data arose from the model postulated.
Hence, we get
\begin{align}
	\bm \varphi^{*} = \arg \max_{\bm \varphi} \log p(\mathcal{D}_N \mid \bm \varphi).
\end{align}
As this is analytically intractable, we propose employing Variational Bayesian  Expectation-Maximization (VB-EM)  \citep{beal_variational_2003} according to which a lower-bound $\mathcal{F}$  to the log-evidence  (called Evidence Lower BOund, ELBO) is constructed with the help of auxiliary densities $h_i(\bm b^{(i)})$ on the unobserved variables $\bm b^{(i)}$: 
\be
\begin{array}{l}
\log p(\mathcal{D}_N \mid \varphi)  \geq   \\
	 \underbrace{\sum_{i=1}^N \E_{h_i(\bm b^{(i)})}\left[ \log \frac{ p(\hat{\bz}^{(i)} \mid \bz(\bs{b}^{(i)}) )  p(\bm b^{(i)}\mid\bm x^{(i)};\bm \varphi)}{h_i(\bm b^{(i)})} \right]}_{=\mathcal{F}(h_{1:N},\bm \varphi)}
	\label{eq:elbo}
\end{array}
\ee
This suggests the following iterative scheme where one alternates between the steps:
\begin{itemize}
	\item \textbf{E-step}: 
	Fix $\bm \varphi$ and maximize $\mathcal{F}$  with respect to $h_i(\bs{b}^{(i)})$. It can be readily shown \citep{bishop2006pattern} that optimality is achieved by the conditional posterior i.e.   
	\begin{align}\label{eq:E_step}
		h_i^{opt}(\bm b^{(i)}) = p(\bm b^{(i)}\mid\mathcal{D}_N, \bm \varphi) \propto p(\hat{\bz}^{(i)} \mid \bm b^{(i)})\nonumber\\ p(\bm b^{(i)} \mid \bm x^{(i)}, \bm \varphi)
	\end{align}
	which makes the inequality in \refeq{eq:elbo} tight. Since the  likelihood is not tractable as it involves a physics-based solver, we have used Markov Chain Monte Carlo (MCMC) to sample from the conditional posterior (see   \refsec{sec:numericalexperiments})
	\item \textbf{M-step}: Given $\{h_i(\bs{b}^{(i)})\}_{i=1}^N$,  maximize  $\mathcal{F}$ with respect to $\bm \varphi$. 
	\begin{align}
	    \bm{\varphi}^{n+1} = \arg\max_{\bm \varphi} \mathcal{F}(h_{1:N}, \bm{\varphi}^n)
	\end{align}

	This requires   derivatives of $\mathcal{F}$ i.e.: 
	\begin{align}
		\frac{\pa \mathcal{F} }{\pa \bm \varphi} = \sum_{i=1}^N \E_{h_i}\left[\frac{\pa \log p(\bm b^{(i)}\mid\bm{ x}^{(i)};\bm \varphi)}{\pa \bm \varphi} \right]
	\end{align}
	Given the MCMC samples $\{ \bs{b}^{(i)}_m\}_{m=1}^M$ from the E-step, these can be approximated as:
	\begin{align}\label{eq:M_step_approx}
		\frac{\pa \mathcal{F}}{\pa \varphi}  \approx \sum_{i=1}^N \frac{1}{M} \sum_{m=1}^M  \frac{\pa \log p(\bm{b}_m^{(i)}\mid\bm x^{(i)};\varphi)}{\pa \varphi} 
	\end{align}
	Due to the Monte Carlo noise in these estimates,  a stochastic gradient ascent algorithm is utilized.  
	In particular, the ADAM optimizer \citep{kingma2014adam} was used  from the \texttt{PyTorch} \citep{paszke2019pytorch} library to capitalize on its auto-differentiation capabilities. 
\end{itemize}

The major elements of the method are summarized in the Algorithm \ref{Alg:model_learning_algo}. We note  here that training complexity grows linearly with the number of training samples $N$ due to the densities $h_i$ associated with each data point (\textbf{for}-loop of Algorithm \ref{Alg:model_learning_algo}) but this can be embarrassingly parallelized.

\begin{algorithm}
	\caption{Data-based model learning}
	\label{Alg:model_learning_algo}
	\begin{algorithmic}[1]
		\State \textbf{Input}: Data $\mathcal{D}_N=\{ \hat{\bx}^{(i)},  \hat{\bz}^{(i)}\}_{i=1}^N$, model form $p(\bs{b}\mid  \bx, ;~\bm \varphi)$, likelihood noise $\bm{\Sigma}_l$, $n=0$
		\State \textbf{Output}: Learned parameter $\bm{\varphi}^*$
		
		\State Initialize the parameters $\bm \varphi$  
		
		\While{ELBO not converged} \\
		{\textbf{Expectation Step (E-step):}}
		\For{$i = 1$ to $N$}
		\State Sample from  the posterior probability $ p(\bm b^{(i)}\mid\mathcal{D}_N, \bm{\varphi}^n)$ using current $\bm{\varphi}^n$ using MCMC \Comment{Eq. \eqref{eq:E_step}}
		\EndFor
		
		{\textbf{Maximization Step (M-step):}}
		\State Monte Carlo gradient estimate \Comment{Eq. \eqref{eq:M_step_approx}}
		\State $\bm{\varphi}^{n+1} = \arg\max_{\bm{\varphi}} \mathcal{F}(h_{1:N}, \bm{\varphi}^n)$
		\State $n\gets n+1$
		\EndWhile
	\end{algorithmic}
\end{algorithm}
 \textbf{Model Predictions}: The VB-EM based model learning scheme discussed above can be carried out in an  \textit{offline phase}. Once the model is learnt, we are interested in the proposed models ability to produce probabilistic predictions  (\textit{online stage}), that reflect the various sources of uncertainty discussed previously. For learnt parameters $\bm{\varphi}^*$, the predictive posterior density $p_{\text{pred}}(\bz\mid \calD, \bm{\varphi}^*)$ on the solution vector $\bz$ of a physical model is as follows:
 \begin{align}\label{eq:predictive_density}
    p_{\text{pred}}(\bz\mid \calD, \bm{\varphi}^*) &= \int p(\bz, \bb\mid \calD, \bm{\varphi}^*) d\bb \nonumber\\
    &=\int p(\bz\mid\bb)~p(\bb\mid \calD, \bm{\varphi}^*)~d\bb \\
    &     \approx \frac{1}{K} \sum_{k=1}^K \bz( \bm b^{(k)})
 \end{align}
 The second of the densities  
 is the conditional (\refeq{eq:prior}) substituted with the learned $\bm{\varphi}^*$ and the first of the densities is simply a Dirac-delta that corresponds to the solution of the physical model, i.e. $\bz(\bm b)$. The intractable integral can be approximated by Monte Carlo using $K$ samples of $\bm b$ drawn from $p(\bb\mid \calD, \bm{\varphi}^*)$.

\subsection{Optimization under uncertainty}\label{sec:optimization}

With the relevant missing links identified as detailed in the previous section, the optimization can be performed on the basis of   \reffig{fig:stochastic graph demonstrator}. 
We  seek to optimize the objective function $\calO$ subject to constraints $\calC = (\calC_1, \ldots, \calC_I)$ that are dependent on uncertain parameters $\bm b$, which in turn are dependent on the design variables $\bx$. 
In this setting, the general parameter-dependent nonlinear constrained optimization problem can be stated as
\begin{align}\label{eq:general_optimization}
	&\min_{\bx} \calO(y_o(\bx,\bs{b})), \nonumber \\
	&\textrm{s.t}~~~ \calC_i(y_{c_i}(\bx,\bs{b}))\leq 0 ,~~ \forall i \in \{1,\ldots,I\} 
\end{align}
where $\bm{x}$ is a $d$ dimensional vector of design variables and $\bm b$ are the model parameter discussed in the previous section. 
It can be observed that the optimization problem is non-trivial because of three main reasons: a) the presence of the constraints (\refsec{sec:constraints_stochastic}) b)  the presence of random variables $\bs{b}$ in  the objective and the constraint(s)  (\refsec{sec:constraints_stochastic}) and, c) non-differentiability of $y_o, y_{c_i}$ and therefore of $\mathcal{O}$ and  $\mathcal{C}_i$. 

\subsubsection{Handling stochasticity and constraints}\label{sec:constraints_stochastic}
Since the solution of the \refeq{eq:general_optimization} depends on the random variables $\bm b$, the objective and constraints are random variables as well and we have to take their random variability into account. We do this by reverting to  a robust optimization problem \citep{ben1999robust,bertsimas2011theory}, with expected values denoted  by $\E[\cdot]$ being the robustness measure to integrate out the uncertainties. In this manner, the optimization problem in \refeq{eq:general_optimization} is reformulated as:
\begin{align}\label{eq:optimization_expectation}
	&\min_{\bx} \E_{\bm b}[\calO(y_o(\bx,\bs{b}))], \nonumber \\
	&\textrm{s.t}~~~ \E_{\bm b}[\calC_i(y_{c_i}(\bx,\bs{b}))]\leq 0 ,~~ \forall i \in \{1,\ldots,I\} 
\end{align}
The expected objective value will yield a design that performs best on average while the reformulated constraints imply   feasibility on average.
 
One can cast this constrained problem to an unconstrained one using penalty-based methods \citep{wang_stochastic_2003, nocedal1999numerical}. In particular we define an augmented objective function $\mathcal{L}$ as follows: 
\begin{align}\label{eq:updated_objective}
	\mathcal{L}(\bx, \bs{b},\bm{\lambda}) =  \calO(y_o(\bx,\bs{b}))+\nonumber \\ \sum_{i}\lambda_i \max\left(\calC_i(y_{c_i}(\bx,\bs{b})),0\right)
\end{align}
where $\lambda_i>0$ is the penalty parameter for the $i^{th}$ constraint. The larger the $\lambda_i$'s are, the more strictly the constraints are enforced. 
Incorporating the augmented objective (\refeq{eq:updated_objective}) in the reformulated optimization problem (\refeq{eq:optimization_expectation}), one can arrive at the following penalized objective:
\be
\E_{\bs{b}} [\mathcal{L}(\bx, \bs{b},\bm{\lambda}) ] = \int \mathcal{L}(\bx, \bs{b},\bm{\lambda}) p(\bm b \mid \bx, \bm{\varphi})d\bm{b},
\ee
leading to the following equivalent, unconstrained optimization problem:
\begin{align}\label{eq:updated_opt_prob}
    \min_{\bx} \E_{\bs{b}} [\mathcal{L}(\bx, \bs{b},\bm{\lambda}) ].
\end{align}
The expectation above is approximated by Monte Carlo which induces noise and necessitates the use of stochastic optimization methods (discussed in detail in the sequel).
Furthermore, we propose to alleviate the dependence on the penalty parameters $\bm \lambda$ by using the sequential unconstrained minimization technique (SUMT) algorithm \citep{fiacco1990nonlinear}, which has been shown to work with non-linear constraints \citep{liuzzi2010sequential}.  
The algorithm considers a strictly increasing sequence $\{\bm \lambda_n\}$ with $\bm \lambda_n \rightarrow \infty$. \citep{fiacco1990nonlinear} proved that when $\bm \lambda_n \rightarrow \infty$, then the sequence of corresponding minima, say $\{\bx^*_n\}$, converges to a global minimizer $\bx^*$ of the original constrained problem. This adaptation of the penalty parameters helps to balance the need to satisfy the constraints with the need to make progress towards the optimal solution. 

\subsubsection{Non-differentiable objective and constraints}\label{sec:non-diff_obj_cons}
We note that the approximation of the objective in \refeq{eq:updated_opt_prob} with Monte Carlo requires multiple runs of the expensive, forward, physics-based models involved, at each iteration of the optimization algorithm. In order to reduce the number of iterations required, especially when the dimension of the design space is higher, derivatives of the objective would be needed. 
In cases where the dimension of the design vector $\bs{x}$ is high, gradient-based methods are necessary. 
In turn, the computation of  derivatives of $\mathcal{L}$ would necessitate derivatives of the outputs of the forward models with respect to the optimization variables $\bm x$.  
The latter are however  unavailable due to the non-differentiability of the forward models. This is a common, restrictive feature of several physics-based simulators which in most cases of engineering practice are implemented in legacy codes that are run as black boxes. This lack of differentiability  has been recognized as a significant roadblock by several researchers   in recent years  \citep{cranmer2020frontier, louppe_adversarial_2019, beaumont2002approximate,marjoram2003markov, agrawal2023probabilistic,shirobokov2020black, lucor2022simple}. In this work, we advocate  Variational Optimization \citep{bird_stochastic_2018,staines2013optimization}, which employs a differentiable bound on the non-differentiable objective.  
In the context of the current problem, we can write:
\begin{align}\label{eq:VO_main_eq}
	\min_{\bx} &\int \left(\mathcal{L}(\bx,\bs{b},\bm{\lambda})\right) p(\bm b \mid \bx, \bm{\varphi}) d\bm{b} \leq \nonumber\\
	&\underbrace{\int \left(\mathcal{L}(\bx,\bs{b},\bm{\lambda})\right) p(\bm b \mid \bx, \bm{\varphi}) q(\bx \mid \bt)d\bm{b}d\bx}_{U(\bt)}
\end{align}
where $q(\boldsymbol{x}\mid \bs{\theta})$ is a  density over the design variables $\bx$ with parameters $\bm{\theta}$. If $\bx^*$ yields the minimum of the objective $\E_{\bm b}[\calL]$, then this can be achieved with a degenerate $q$ that collapses to a Dirac-delta, i.e. if $q(\boldsymbol{x}\mid \bs{\theta})=\delta(\bx-\bx^*)$. For all other densities $q$ or parameters $\bt$, the inequality above would in general be strict. 
Hence and instead of minimizing $\E_{\bm b}[\calL]$ with respect to $\bm{x}$, we can minimize the upper bound $U$ with respect to $\bm{\theta}$. Under mild restrictions outlined by \citep{staines_variational_2012}, the bound $U(\boldsymbol{ \theta})$ is differential w.r.t $\bm{\theta}$ and using the log-likelihood trick its gradient can be rewritten as \citep{williams1992simple}:
\begin{align}\label{eq:grad_estimator}
	\nabla_{\bt} &U(\bt) = 
	\E_{\bx, \bm{b}}\left[\nabla_{\bt} \log q(\bx \mid \bt) \mathcal{L}(\bx,\bs{b},\bm{\lambda})\right] \\
	& \approx \frac{1}{S} \sum_{i=1}^{S} \mathcal{L}(\bx_i,\bs{b}_i,\bm{\lambda}) \frac{\partial}{\partial \bt} \log q\left(\bx_i \mid \theta\right)
\end{align}
%
Both terms in the integrand can be readily evaluated which opens the door for a Monte Carlo approximation of the aforementioned expression by drawing  samples $\bx_i$ from $q(\bx \mid \bt)$ and  subsequently $\bm b_i$  from $p(\bm b \mid  \bx_i, \varphi^*)$.

\subsubsection{Implementation considerations}
While the Monte Carlo estimation  of the gradient of the new objective $U(\bt)$ also requires running the expensive, forward models multiple times, it can be embarrassingly parallelized.

Obviously, convergence is impeded by the unavoidable Monte Carlo errors in the aforementioned estimates. In order to reduce them, we advocate the  use of the baseline estimator proposed in \citep{kool_buy_2022} which is based on  the following expression:
\begin{align}
	\frac{\partial U}{\partial \theta} \approx &\frac{1}{S-1} \sum_{i=1}^{S}  \frac{\partial}{\partial \bt} \log q\left(\bx_i \mid \bt\right) \nonumber \\
	&\left(\mathcal{L}(\bx_i,\bs{b}_i,\bm{\lambda}) - \frac{1}{S}\sum_{j=1}^S \mathcal{L}(\bx_j,\bs{b}_j,\bm{\lambda}) \right) \label{eq:baseline_trick}
\end{align}
The estimator above is also unbiased as the one in \refeq{eq:grad_estimator}, it 
 does not imply any additional cost beyond the  $S$ samples and in addition exhibits lower variance as shown in \citep{kool_buy_2022}.

\begin{figure*}[!htbp]
    \centering
  \includegraphics[width=0.8\textwidth]{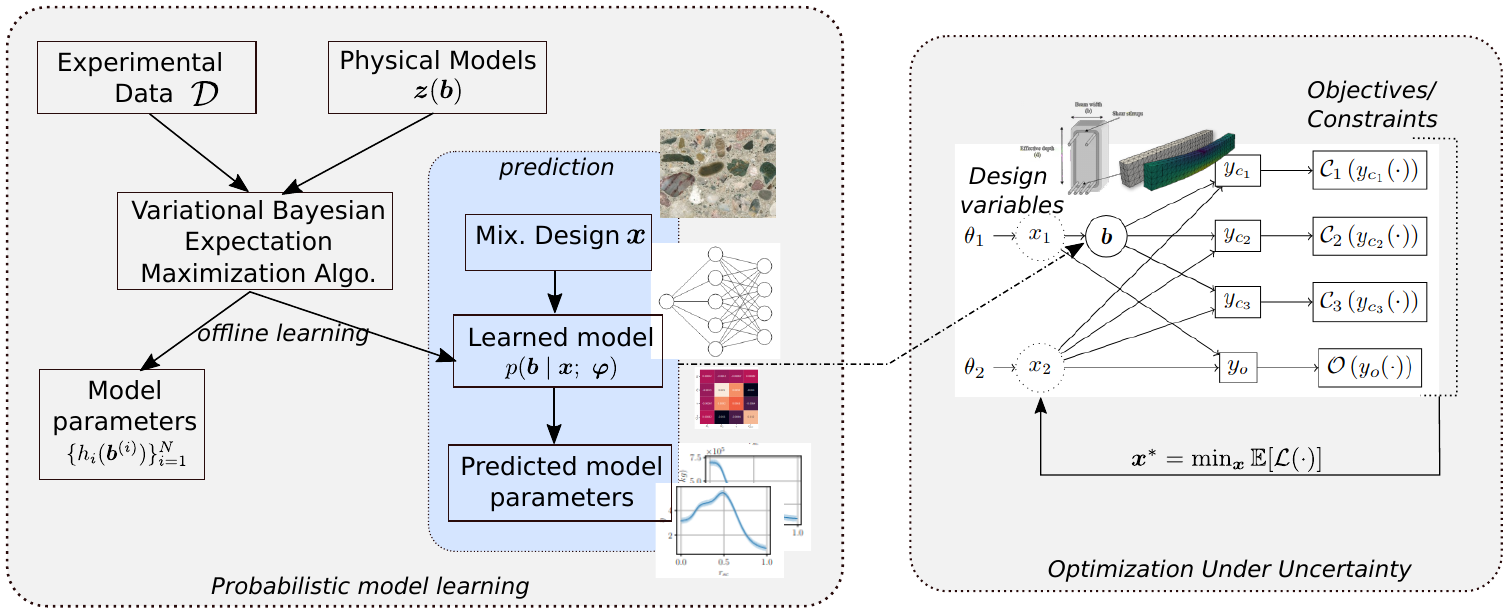}
  \caption{A schematic of the probabilistic model learning (left block) and the optimization under uncertainty (right block). The left block illustrates how the information from experimental data and physical models are fused together to learn the missing probabilistic link. This learned probabilistic link subsequently becomes a linchpin in predictive scenarios, particularly in downstream optimization tasks. The right block illustrates querying the learned probabilistic model to complete the missing link and interfacing the workflow describing the design variables, physical models and the KPIs. Subsequently, this integrated approach facilitates the execution of optimization under uncertainty as per the proposed methodology}
  \label{fig:schematic_summary}
\end{figure*}

To efficiently compute the gradient estimators, we make use of the auto-differentiation capabilities of modern machine learning libraries. In the present study, \texttt{PyTorch} \citep{paszke2019pytorch} was utilized. For the stochastic gradient descent, the ADAM optimizer was used \citep{kingma2014adam}. In the present study, $q(\boldsymbol{x}\mid \bs{\theta})$ was a  Gaussian distribution with parameters $\bs{\theta} = \{ \bm{\mu} , \bm{\Sigma}\}$ representing mean and diagonal covariance, respectively. We say we have arrived at an optimal $\bx^*$ when the $q$ almost degenerates to a Dirac-delta, or colloquially, when the variance of $q$ has converged and is considerably small. For completeness, the algorithm for the proposed optimization scheme is given in Algorithm \ref{Alg:optimization_algo}. A schematic overview of the methods discussed in \refsec{sec:model_learning} and the \refsec{sec:optimization} is presented in \reffig{fig:schematic_summary}.
\begin{algorithm}[h]
	\caption{Black-box stochastic constrained optimization}
	\label{Alg:optimization_algo}
	\begin{algorithmic}[1]
		\State \textbf{Input}: distribution $q(\bm{x} \mid \bt)$ for the design variable $\bx$, $n = 0$, learning rate $\eta$
		\State $\bt_0^0, \bm{\lambda}_1 \gets$ choose starting point
		\For{$k = 1, 2, \ldots$}
		\While{$\bt_k$ not converged}
		\State Sample design variables and model parameters
  $\bm{x}_i \sim q(\bm{x} \mid \bt_k^n), \quad \bm{b}_i \sim p(\bm b \mid \bx_i, \bm{\varphi})$
		\State Farm the workflow with physics-based solvers for the samples in different machines and compute updated objective $\mathcal{L}(\bx_i,\bs{b}_i,\bm{\lambda}_k)$ for each of them     \Comment{Eq. \eqref{eq:updated_objective}}
		\State Compute baseline
		\State  
            Monte-Carlo gradient estimate $\nabla_{\bt} U$ \Comment{Eq. \eqref{eq:baseline_trick}}
		\State $\bt_k^{n+1} \gets \bt_k^n + \eta \nabla_{\bt} U$ \Comment{Stochastic Gradient Descent}
		\State $n \gets n + 1$
		\EndWhile
		\If{$\|\bt_{k}^n - \bt_{k-1}^n\| \leq \varepsilon$} \Comment{Convergence condition}
		\State \textbf{break}
		\EndIf
		\State $\bm{\lambda}_{k+1} \gets \bm{\lambda}_{k}$ \Comment{Update penalty parameter} 
		\State  $\bt_{k+1}^{0} \gets \bt_{k}^n$ \Comment{Update starting parameter}
		\EndFor
		\State \Return{$\bt_{k}$}
	\end{algorithmic}
\end{algorithm}

\section{Numerical Experiments}\label{sec:numericalexperiments}
This section presents the results of the data-based model learning (\refsec{sec:model learning results}) and optimization (\refsec{sec:optimization results}) methodological frameworks for the coupled workflow describing the concrete mix-design and structural performance as discussed in \reffig{fig:snakemake_workflow}. 

\subsection{Model-learning results}\label{sec:model learning results}
We report on the results obtained upon application of the method presented in \refsec{sec:model_learning} on the hydration and homogenization models (Table \ref{table:model_learning}) along with the experimental data (\refsec{sec:exp_data}). 

\textbf{Implementation details:} We note that for the likelihood model (\refeq{eq:likelihood}) corresponding to the hydration model, the observed output $\hat{\bz}$ in the observed data $\calD =\{\hat{\bx}^{(i)},  \hat{\bz}^{(i)}\}_{i=1}^N$  is the cumulative heat $\hat{\bm{Q}}$ and the covariance matrix $\Sigma_{\ell} = 3^2\bm{I}_{dim(\hat{\bm{Q})}\times dim(\hat{\bm{Q})}}$. 
The particular choice was made to account for the cumulative heat sensor noise of $\pm 4.5 \unit{J/g}$ as reported in \citep{gruyaert2011}. For the homogenization model, the 
$\hat{\bz}= [E_c,f_c]^T$ where $E_c$ is the Young's Modulus and $f_{c}$ is the compressive strength of the concrete. The covariance matrix $\Sigma_{\ell} = \text{diag}(4\times10^{18},2\times 10^{12})\unit{Pa^2}$. For both of the above, $\hat{\bx}$ in the observed data $\calD$ is the slag-binder mass ratio $r_{sb}$. 

For both cases, a fully-connected neural network is used to parameterize the mean of the conditional of model parameters $\bb$ (\refeq{eq:prior}). The optimum number of hidden layers and nodes per layer was determined to be 1 and 30 respectively. The Tanh was chosen as activation function for all layers. The $L_2$ weight regularization was employed to prevent over-fitting. We employed a learning rate of $10^{-2}$ for all the results reported here.

Owing to the intractability of the conditional posterior given in \refeq{eq:E_step}, we approximate it with MCMC, in particular we used the DRAM (Delayed Rejection Adaptive Metropolis) \citep{2020arXiv200809589S, haario2006dram}. The specific selection was motivated by two primary considerations. Firstly, a gradient-free sampling strategy is imperative due to the absence of gradients in the physics-based models employed in this context. Secondly, to introduce automation to the tuning of free parameters in the Markov Chain Monte Carlo (MCMC) methods, ensuring a streamlined and efficient convergence process. In the DRAM sampler, we bound the target acceptance rate to be between $0.1$ to $0.3$. 

\begin{figure}[!htpb]
	\centering
	\includegraphics[width=0.5\textwidth]{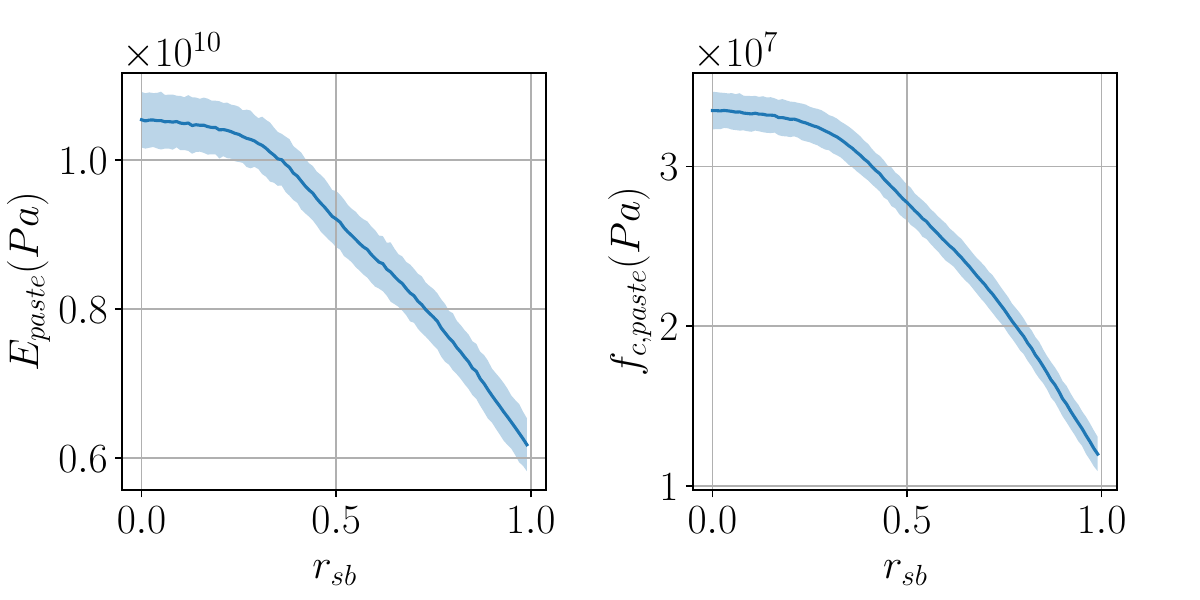}
	\caption{Learned probabilistic relation between the homogenization model parameters  and the slag-binder ratio $r_{sb}$. The solid line denotes the mean and the shaded $\pm 2 \times$ standard deviation}
	\label{fig:learned_model_homo}
\end{figure}

\begin{figure}[!htpb]
	\centering
	\includegraphics[width=0.5\textwidth]{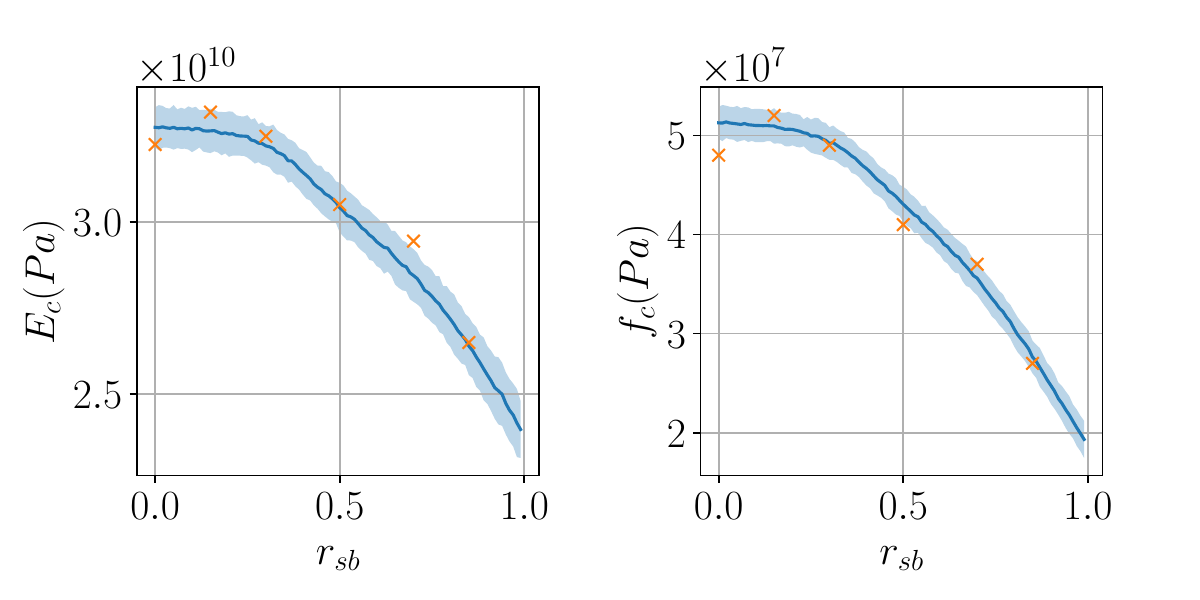}
	\caption{Predictive performance of the learned model corresponding to the homogenization process. The Solid line is the predictive mean, and the shaded area is $\pm 2 \times$ standard deviation. The crosses correspond to the noisy observed data}
	\label{fig:prediction_homo_model}
\end{figure}

\textbf{Results:} \reffig{fig:learned_model_homo} shows the learned probabilistic relation between the latent model parameters of the homogenization model and the slag-binder mass ratio $r_{sb}$. Out of the six available noisy datasets (\refsec{sec:exp_data}), five were used for training and the dataset corresponding to $r_{sb} = 0.5$ was used for testing. We access the predictive capabilities of the learned model by propagating the uncertainties forward via the homogenization model and analyzing the predictive density $p_{pred}$ (\reffig{eq:predictive_density}) as illustrated in \reffig{fig:prediction_homo_model}. We observe that the mechanical properties of concrete obtained by the homogenization model with learned probabilistic model predictions as the input, envelops the ground truth.

\begin{figure}[!htpb]
	\centering
	\includegraphics[width=0.5\textwidth]{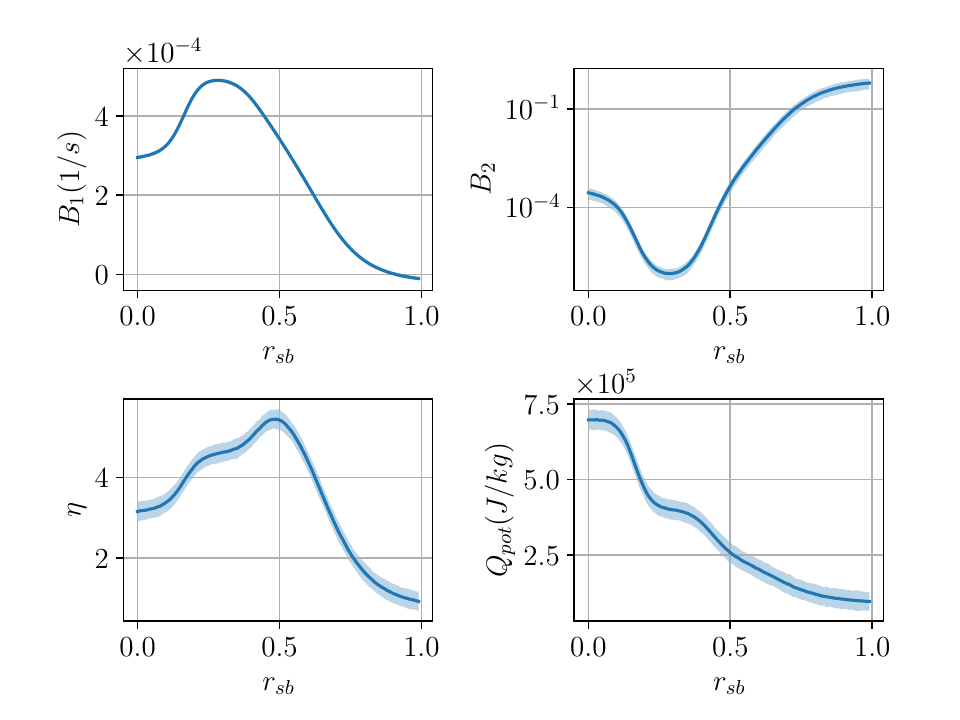}
	\caption{Learned probabilistic relation between the hydration model parameters and the slag-binder mass ratio $r_{sb}$. Solid line denotes the mean and the shaded $\pm 2 \times$ standard deviation}
	\label{fig:learned_model_hydration}
\end{figure}

\begin{figure}[!htpb]
	\centering
	\includegraphics[width=0.4\textwidth]{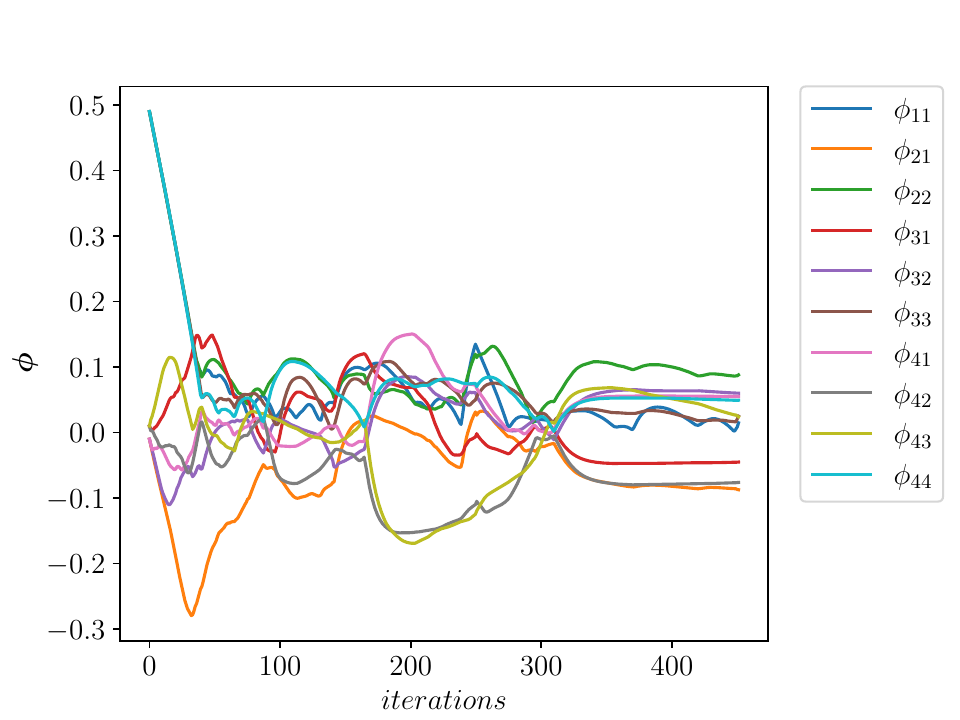}
	\caption{Evolution of the entries $\phi_{ij}$ of the lower-triangular matrix $\bm{L}$ of the covariance matrix (\refeq{eq:prior}) with respect to EM iterations}
	\label{fig:cov_evolution}
\end{figure}

\begin{figure}[!htpb]
	\centering
	\includegraphics[width=0.4\textwidth]{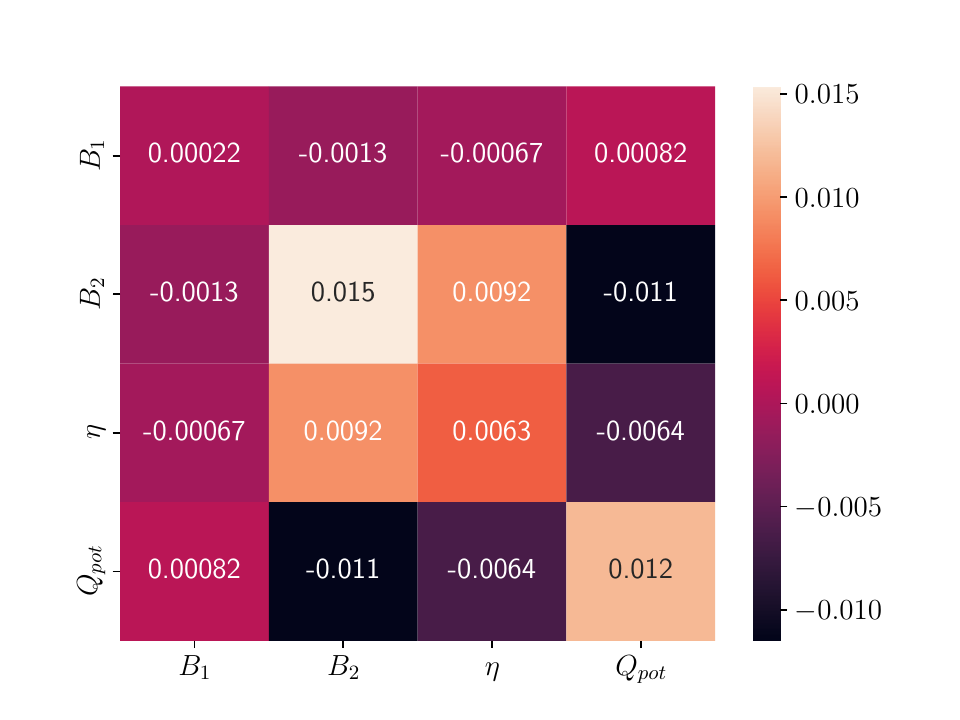}
	\caption{Heat map of the converged value of the covariance matrix $\bs{LL}^T$ of the probabilistic model corresponding to concrete hydration}
	\label{fig:cov_heat_map}
\end{figure}

\begin{figure}[!htpb]
	\centering
	\includegraphics[width=0.45\textwidth]{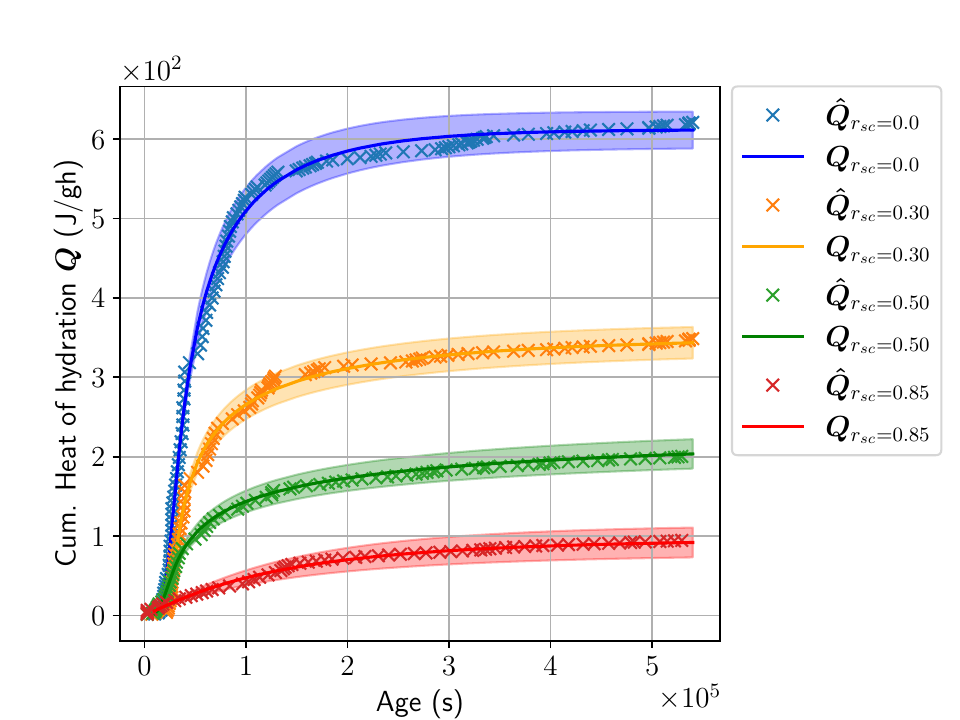}
	\caption{Predictive performance of the learned model corresponding to the hydration process. The solid line is the predictive mean and the shaded area is $\pm2\times$ standard deviation. The crosses correspond to the noisy observed data}
	\label{fig:prediction_hyd_model}
\end{figure}

Similarly, for the hydration model, \reffig{fig:learned_model_hydration} shows the learned probabilistic relation between the latent model parameters $(B_1,B_2,\eta,Q_{pot})$ and the ratio $r_{sb}$. Out of the four available noisy datasets (\refsec{sec:exp_data}) for $T = 20^{\circ} \unit{C}$, three were used for training and the dataset corresponding to $r_{sb} = 0.5$ was used for testing. The value of $\activE$ was taken from \citep{gruyaert2011}. \reffig{fig:prediction_hyd_model} compares the experimental heat of hydration for different $r_{sb}$  with the probabilistic predictions made using the learned probabilistic model as an input to the hydration model. We observe that the predictions entirely envelop the ground truth data, while accounting for the aleatoric noise present in the experimental data. \reffig{fig:cov_evolution} shows the evolution of the entries of the covariance matrix of the conditional on the hydration model latent parameters $p(\bs{b}\mid \bx; ~\bm{\varphi})$. It serves as an indicator for the convergence of the EM algorithm. The converged value of the covariance matrix is given by \reffig{fig:cov_heat_map}. It confirms the intricate correlation among the hydration model parameters, also reported in \reffig{fig:heatrelease}. This is a general challenge with most physical models that are often overparameterized (at least for a given data set) leading to multiple configurations of parameters with similar likelihood (see  \reffig{fig:parameterEvolution}).

At this point, it is crucial to (re)state that the training is performed using {\em indirect},   noisy data. It is encouraging to note that the learned models are able to account for the aleatoric uncertainty arising from the noise in the observed data and the epistemic uncertainty due to the finite amount of training data. The probabilistic model is able to learn relationships which were otherwise unavailable in literature, with the aid of physical models and (noisy) data.

\subsection{Optimization results}\label{sec:optimization results}

\begin{figure*}[!htpb]
    \centering
    \subfigure[]{\includegraphics[width=0.31\textwidth]{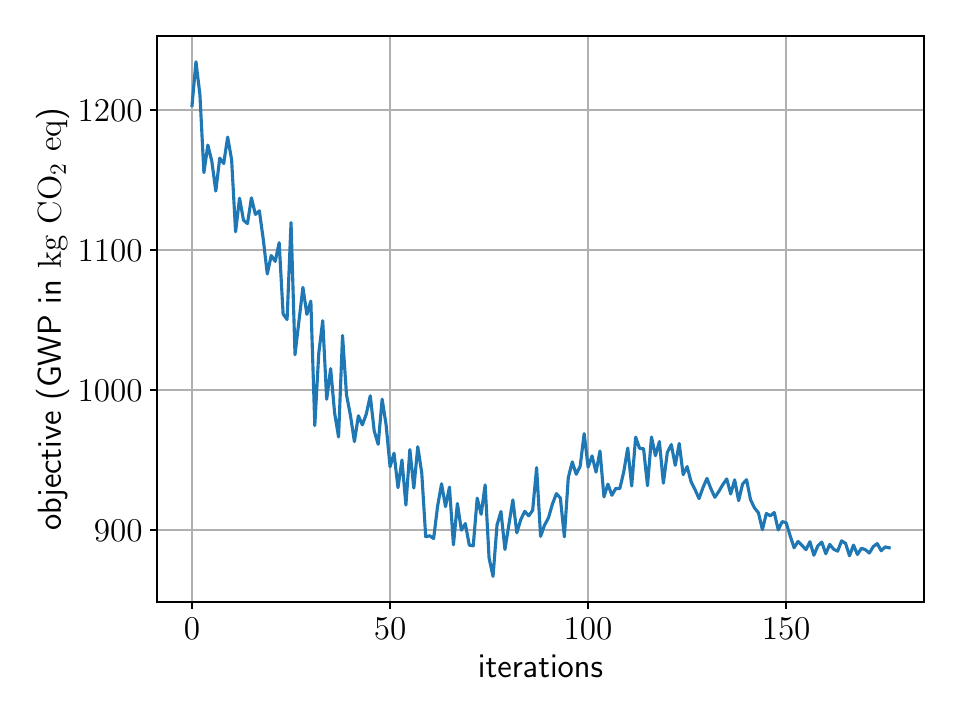}\label{fig:objective_evo}}
    \subfigure[]{\includegraphics[width=0.31\textwidth]{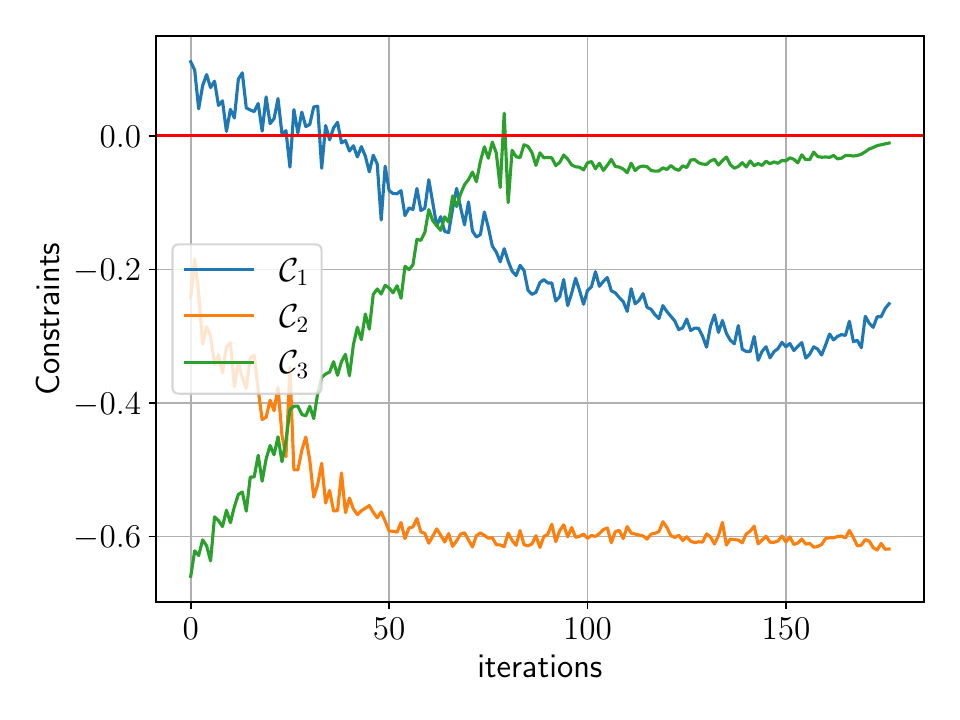}\label{fig:constraint_evo}} 
    \subfigure[]{\includegraphics[width=0.31\textwidth]{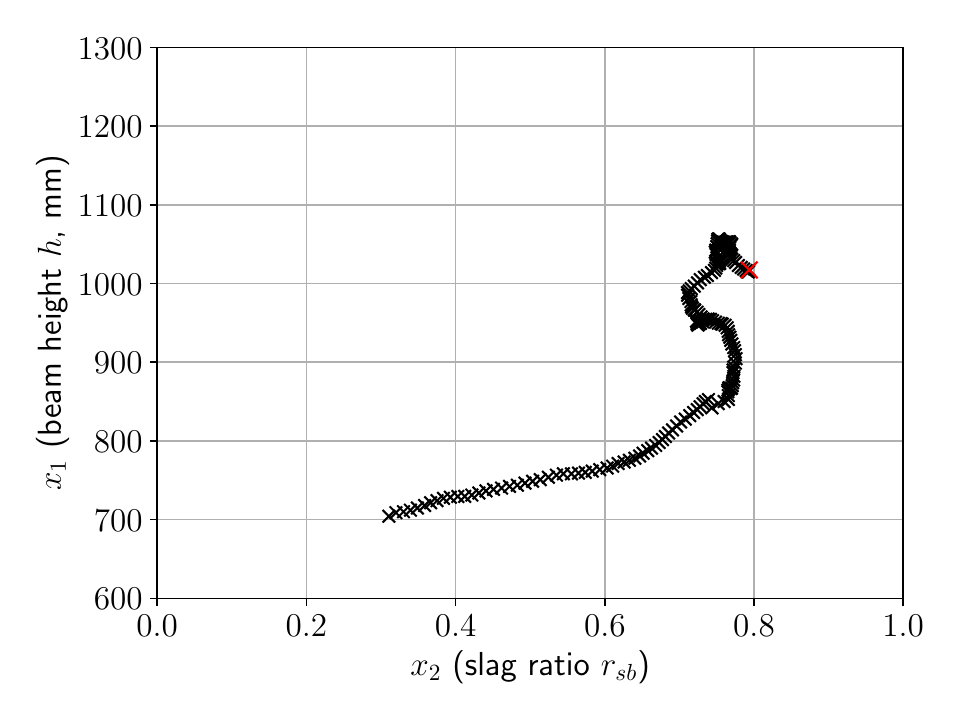}\label{fig:deisgn_var_evo}}
    \caption{(a) Evolution of the expected objective $\E_{\bb}[\calO]$ vs the number of iterations. The objective is the GWP of the beam. (b) Evolution of the expected constraints $\E_{\bb}[\calC_i]$ (which should all be negative) vs the number of iterations. $\calC_1$ represents the beam design constraint, $\calC_2$ the temperature constraint and $\calC_3$ gives the demoulding time constraint. 
    (c) Trajectory of the design variables (slag-binder mass ratio $r_{sb}$ and the beam height $h$). The red cross represents the optimal  value identified upon convergence}
    \label{fig:optimization_results}
\end{figure*}

\begin{figure}[!htpb]
	\centering
	\includegraphics[width=0.3\textwidth]{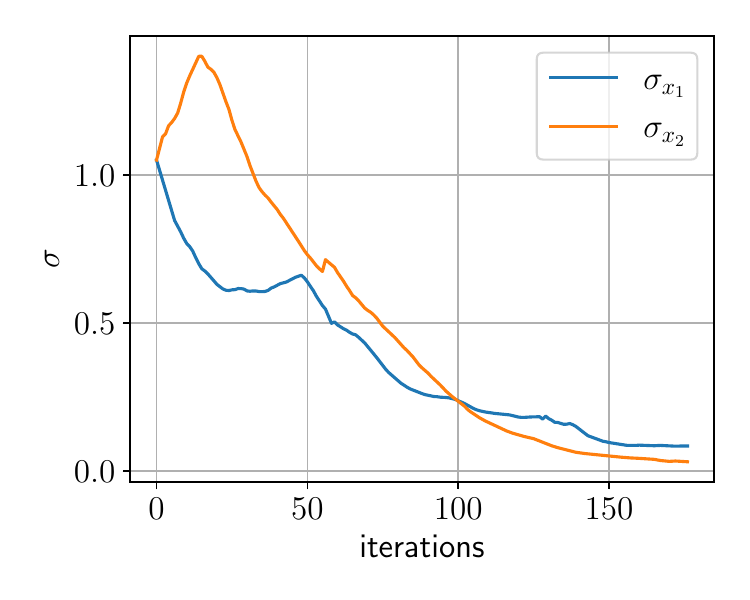}
	\caption{Evolution of the standard deviations $\sigma$ of the design variables to highlight convergence of the optimization process. We note that the design variables are transformed and scaled in the optimization procedure}
	\label{fig:opt_convergance}
\end{figure}
With the learned probabilistic links as discussed in the previous section, we overcame the issue of forward and backward information flow bottleneck in the workflow connecting the design variables and KPIs relevant for constraints/objetcive (as discussed in \refsec{sec:calibration}). In this section, we report on the results obtained by performing optimization under uncertainty as discussed in \refsec{sec:optimization} for the performance-based concrete design workflow. The design variables, objectives, and the constraints are detailed in Table \ref{table:stoc_graph_variables}. For the Temperature constraint, we choose $T_{limit} = 60 \degree C$ and for the demoulding time constraint, we choose 10 hours. 
To improve the numerical stability, We scale the variables, constraints, and objectives to make them of the order 1. 
To demonstrate the optimization scheme proposed, a simply supported beam is used as discussed in \refsec{sec:example_problem}, with parameters given in Table \ref{tab:beamdesigninput}. Colloquially, we aim to find the value(s) of slag-binder mass ratio $r_{sb}$ and beam height $h$ that minimize the objective, on average, while satisfying, on average, the aforementioned constraints.

As discussed, the workflow for estimating the gradient is embarrassingly parallelizable.
Hence, for each value considered in the design space, we call the ensemble of the workflows in parallel machines and  collect the results. 
For the subsequent illustrations, a step size of $0.05$ is utilized in the ADAM optimizer and $S=100$ was the number of samples for gradient estimation.  We set $\lambda_i = 1 ~\forall i \in \{1,\ldots,I\} $ as the starting value. \reffig{fig:optimization_results} shows the optimization results. In the design space, we start from values of design variables that violate the beam design constraints $\calC_1$ as evident from \reffig{fig:constraint_evo}. 
This activates the corresponding penalty term in the augmented objective (\refeq{eq:updated_objective}), thus driving the design variables to satisfy the constraint (around iteration 40). 
Physically, this  implies that the beam is not able to withstand the applied load for the given slag-binder ratio, beam height and other material parameters (which are kept constant in the optimization procedure). As  a result, the optimizer suggests to increase the beam height $h$ in order to satisfy the constraint, while simultaneously increasing the slag-binder mass ratio $r_{sb}$ also, owing to the influence of the GWP objective ( see \reffig{fig:deisgn_var_evo}). 
As it can be seen in \reffig{fig:objective_evo}, this leads to a reduction of the GWP because with the increase of the slag ratio, the Portland cement content which is mainly responsible for the $CO_2$ emission is ultimately reduced. In theory, the optimum value of the slag-binder mass ratio $r_{sb}$ approaches one (meaning only slag in the mix)
 if only the GWP objective with no constraints were to be used in the optimization. But the demoulding time constraint $\calC_3$ penalizes the objective to limit the slag-binder $r_{sb}$ ratio to be around $0.8$ (see \reffig{fig:deisgn_var_evo}), since the evolution of mechanical properties is both much faster for Portland cement than for slag and at the same time the absolute values for strength and Young's modulus are higher. 
  This is also evident in \reffig{fig:constraint_evo}, when around iteration $80$, the constraint violation line is crossed thus activating the penalty from $\calC_3$. This also stops the close to linear decent of the GWP objective. In the present illustrations, a  value of $10$ hours is chosen as the demoulding time to demonstrate the procedure. In real-world settings, a manufacturer would be inclined to remove the  formwork earlier so that it can be reused.  But the lower the requirement of the demoulding time, the higher the ratio of cement content required in the mix, leading to an increased hydration heat which in effect accelerates the reaction.

The oscillations in the objective and the constraints as seen in \reffig{fig:objective_evo} and \reffig{fig:constraint_evo} are due to the  Monte Carlo noise in the gradient estimation. As per \refeq{eq:VO_main_eq}, the design variables are treated as random variables following  a normal distribution. As discussed in Algorithm \ref{Alg:optimization_algo},  the optimization procedure is assumed to converge when the    standard deviations $\sigma$ of the normal distribution attain small values (\reffig{fig:opt_convergance}) i.e. when the normal degenerates to a Dirac-delta.  
The $\sigma$ values stabilizing to relatively small values points towards the convergence of the algorithm.
 
The performance increase (in terms of GWP) is difficult to evaluate in the current setting. This is due to the fact that the constraint $\mathcal{C}_1$ is not fulfilled for the initial value of the design variables chosen for the optimization. It is to be highlighted that this is actually an advantage of the method - the user can start with a reasonable design that still violates the constraints. In order to make a reasonable comparison, a design using only Portland cement (i.e., $r_{sb}$ = 0) with only the loadbearing capacity as a constraint (beam design constraint $\mathcal{C}_1$) and the free parameter being the height of the beam was computed. This minimum height was found to be 77.5cm with a corresponding GWP of the beam of \mbox{1455 $\unit{kg CO_2 eq}$}. Note that this design does not fulfill the temperature constraint $\mathcal{C}_2$ with a maximum temperature of \mbox{81$\unit{\degree C}$}.
Another option for comparison is the first iteration number in the optimization procedure that fulfills all the constraints in expectation, which is the iteration number $30$ with a GWP of \mbox{1050 $\unit{kg CO_2 eq}$}.
In the subsequent iteration steps, this is further reduced to \mbox{~900 $\unit{kg CO_2 eq}$} for the optimum value of the design variables obtained in the present study.
This reduction in GWP is achieved by increasing the height of the beam 
to 100$\unit{cm}$ while replacing portland cement with blast furnace slag so that the mass fraction of slag-binder $r_{sb}$ is 0.8. 
The addition of slag in the mixture decreases the strength of the material 
as illustrated in \reffig{fig:prediction_homo_model}, while at the same time this decrease is compensated by an increased height. It is also informative to study the evolution of the (expected) constraints shown in \reffig{fig:constraint_evo}. One observes that $\mathcal{C}_3$ (green line) associated with the demoulding time is the most critical. Thus, in the current example, the GWP could be decreased even further when the time of demoulding is extended (depending on the production process of removing the formwork).



\section{Conclusion and Outlook}\label{sec:results}
We introduced a systematic  design approach for precast concrete industry in the pursuit of sustainable construction practices. It makes use of a holistic optimization framework which combines concrete mixture design with the structural simulation of the precast concrete element within an automated workflow. In this manner various objectives and constraints such as the environmental impact of the concrete element or its  structural efficiency, can be considered.

The proposed holistic approach  is  demonstrated on  a specific design problem, but should serve as template that can be readily adapted to other design problems. 
The advocated black-box stochastic optimization procedure is able to deal with  the challenges presented by general workflows, such as the  presence  of black-box models without  derivatives, the effect of uncertainties and of  non-linear constraints.
Furthermore, to complete the forward and backward information flow that is essential in the optimization procedure, a method to learn missing (probabilistic) links between the concrete mix design variables and model parameters from experimental data is presented. We note that to the best of our knowledge,  such a link is not available in the literature. 

We demonstrated on the precast concrete element the integration of material and structural design in a joint workflow and showcased that this has the potential to decrease the objective, i.e.  the global warming potential. For the structural design, semi-analytical models based on the Eurocode are used, whereas the manufacturing process is simulated using a complex FE-model. This illustrates the ability of the proposed procedure to combine multiple simulation tools of varying complexity, accounting for different parts of the life cycle. Hence, extending this in order to include e.g. additional load configurations, materials or life cycle models, is straightforward. The present approach to treat the design process as a workflow, learning the missing links from data/models and finally using this workflow in a global optimization is transferable to several other materials, structural and mechanical problems.
Such extensions could readily  include more complex design processes with an increased number of parameters and constraints (the latter due to multiple load configurations or limit states in a real structure). Furthermore, this procedure could be applied to problems involving a complete structure (e.g. bridge, building) instead of a single element and potentially entailing advanced modeling features that include multiscale models to link material composition to material properties or improve the computation of the global warming potential using a complete life cycle analysis.

\label{sec:conclusion}



\subsection*{Acknowledgments}
The authors gratefully acknowledge the support of VDI Technologiezentrum Gmbh and the Federal Ministry for Education and Research (BMBF) within the collaborative project  "Lebenszyklus von Beton - Ontologieentwicklung für die Prozesskette der Betonherstellung". AA and PSK received support through the subproject "Probabilistic machine learning for the calibration and validation of physics-based models of concrete" (FKZ-13XP5125B). E.T and J.F.U received support through the subproject  "Ontologien und Workflows für Prozesse, Modelle und Abbildung in Datenbanken" (FKZ-13XP5125B). We also acknowledge the support of David Alos Shepherd (Building Materials Technology, Karlsruhe Institute of Technology) for digitizing the experimental data.
\subsection*{Statements and Declarations}

\textbf{Competing interests} The authors declare that they have no known
competing financial interest or personal relationships that could have
appeared to influence the work reported in this paper.\\

\noindent \textbf{Replication of results} The paper provides a enough description of the proposed method so that the results can be replicated. If need be, the code and data of this study are available on
request from the corresponding author.


\bibliography{optimization_paper_bibliography_BAM,optimization_paper_bibliography_TUM}
\begin{appendices}
\section{Homogenization}\label{appendix:hom}
\subsection{Approximation of elastic properties}\label{ssec:mt_elastic}
The chosen method to homogenize the elastic, isotropic properties $\eMod$ and $\poission$ is the Mori-Tanaka homogenization scheme, \cite{mor_1973_asi}.
It is a well-established, analytical homogenization method.
The formulation uses bulk and shear moduli $\bulkMod$ and $\shearMod$.
They are related to $\eMod$ and $\poission$ as $\bulkMod = \frac{\eMod}{3(1-2\poission)}$ and $\shearMod = \frac{\eMod}{2(1+\poission)}$.
The used Mori-Tanaka method assumes spherical inclusions in an infinite matrix and considers the interactions of multiple inclusions.
The applied formulations follow the notation published in 
\cite{nee_2012_ammf}, where this method is applied to successfully model the effective concrete stiffness for multiple types of inclusions.
The general idea of this analytical homogenization procedure is to describe the overall stiffness of a body $\body$, based on the properties of the individual phases, i.e. the matrix and the inclusions.
Each of the $n$ phases is denoted by the index $\phaseIndex$, where $\phaseIndex = 0$ is defined as the matrix phase.
The volume fraction of each phase is defined as
\begin{align}
	\volFracPhase = \frac{\left\| \bodyPhase \right\|}{\left\| \body \right\|} \quad  \text{for}~ \phaseIndex = 0, ..., n.
\end{align}
The inclusions are assumed to be spheres, defined by their radius $\radiusPhase$.
The elastic properties of each homogeneous and isotropic phase is given by the material stiffness matrix $\bL^{(\phaseIndex)}$, here written in terms of the bulk and shear moduli $\bulkMod$ and $\shearMod$,
\begin{align}
	\matStiffPhase= 3 \bulkModPhase \orthProjV + 2 \shearModPhase \orthProjD  \quad \text{for}~ \phaseIndex = 0, ..., n, \label{eq:Lr}
\end{align}
where $\orthProjV$ and $\orthProjD$ are the orthogonal projections of the volumetric and deviatoric components.\\
The method assumes that the micro-heterogeneous body $\body$ is subjected to a macroscale strain $\strain$.
It is considered that for each phase a concentration factor $\concentrationPhase$ can be defined such that
\begin{align}
	\strainPhase = \concentrationPhase\strain \quad  \text{for}~ \phaseIndex = 0, ..., n, \label{eq:strainaverage}
\end{align}
which computes the average strain $\strainPhase$ within a phase, based on the overall strains.
This can then be used to compute the effective stiffness matrix $\matStiffEff$ as a volumetric sum over the constituents weighted by the corresponding concentration factor 
\begin{align}
	\matStiffEff = \sum_{\phaseIndex=0}^{n} \volFracPhase \matStiffPhase\concentrationPhase \quad  \text{for}~ \phaseIndex = 0, ..., n.\label{eq:Leff}
\end{align}

The concentration factors $\concentrationPhase$,
\begin{align}
	\concentrationZero &= \left( \volFracZero\bI + \sum^{n}_{\phaseIndex=1} \volFracPhase \dilConcentrationPhase\right)^{-1}\label{eq:A0}\\
	\concentrationPhase &= \dilConcentrationPhase\concentrationZero\quad  \text{for}~ \phaseIndex = 1, ..., n,
\end{align}
are based on the dilute concentration factors $\dilConcentrationPhase$, which need to be obtained first.
The dilute concentration factors are based on the assumption that each inclusion is subjected to the average strain in the matrix $\strainZero$, therefore
\begin{align}
	\strainPhase = \dilConcentrationPhase\strainZero \quad  \text{for}~ \phaseIndex = 1, ..., n. 
\end{align}
The dilute concentration factors neglect the interaction among phases and are only defined for the inclusion phases $\phaseIndex = 1,...,n$.
The applied formulation uses an additive volumetric-deviatoric split, where
\begin{align}
	\dilConcentrationPhase = \dilConcentrationVPhase\Ivol +  \dilConcentrationDPhase \Idev \quad  \text{for}~ \phaseIndex = 1, ..., n, \text{\quad with}
\end{align}
\begin{align}
	\dilConcentrationVPhase = \dfrac{\bulkModZero}{\bulkModZero + \auxAlphaZero(\bulkModPhase - \bulkModZero)}, \\
	\dilConcentrationDPhase = \dfrac{\shearModZero}{\shearModZero + \auxBetaZero(\shearModPhase - \shearModZero)}.
\end{align}
The auxiliary factors follow from the Eshelby solution as
\begin{align}
	\auxAlphaZero = \frac{1 + \poissionZero}{3(1+ \poissionZero)} \quad\text{and}\quad 
	\auxBetaZero = \frac{2(4 - 5\poissionZero)}{15(1 - \poissionZero)}
\end{align}
where  $\poissionZero$ refers to the Poission's ratio of the matrix phase.
The effective bulk and shear modului can be computed based on a sum over the phases
\begin{align}
\bulkModEff = \dfrac{\volFracZero\bulkModZero + \sum^{n}_{\phaseIndex=1} \volFracPhase \bulkModPhase \dilConcentrationVPhase}{\volFracZero + \sum^{n}_{\phaseIndex=1} \volFracPhase \dilConcentrationVPhase},\label{eq:keff} \\
\shearModEff = \dfrac{\volFracZero\shearModZero + \sum^{n}_{\phaseIndex=1} \volFracPhase \shearModPhase \dilConcentrationDPhase}{\volFracZero + \sum^{n}_{\phaseIndex=1} \volFracPhase \dilConcentrationDPhase}.\label{eq:geff}
\end{align}
Based on the concept of \refeq{eq:strainaverage}, with the formulations \refeq{eq:Lr}, \refeq{eq:Leff} and \refeq{eq:A0}, the average matrix stress is defined as 
\begin{align}
\stressZero = \matStiffZero\concentrationZero {\matStiffEff}^{-1}\stress. \label{eq:matrixstress}
\end{align}
\subsubsection{Approximation of compressive strength}\label{ssec:compressivestrength}
The estimation of the concrete compressive strength $\fc$ follows the ideas of \cite{nev_2018_mcam}.
The procedure here is taken from the code provided in the link in \cite{nee_2012_ammf}.
The assumption is that a failure in the cement paste will cause the concrete to crack.
The approach is based on two main assumptions.
First, the Mori-Tanaka method is used to estimate the average stress within the matrix material $\stressMatrix$. 
The formulation is given in \refeq{eq:matrixstress}.
Second, the von Mises failure criterion of the average matrix stress is used to estimate the uniaxial compressive strength
\begin{align}
	{\fc} = \sqrt{3 \Jtwo},  \label{eq:vonMises}
\end{align}
with $\Jtwo(\stress) = \frac{1}{2} \stressD:\stressD$ and $\stressD = \stress - \frac{1}{3} trace(\stress) \bI$.
It is achieved by finding a uniaxial macroscopic stress $\stress = \begin{bmatrix} -\fcEff & 0 & 0 &0&0&0 \end{bmatrix}\TP$, which exactly fulfills the von Mises failure criterion \refeq{eq:vonMises} for the average stress within the matrix $\stressMatrix$.
The procedure here is taken from the code provided in the link in \cite{nee_2012_ammf}.
First, a $\JtwoTest$ is computed for a uniaxial test stress $\stressTest = \begin{bmatrix} \forceTest & 0 & 0 &0&0&0 \end{bmatrix}\TP$. 
Then the matrix stress $\stressMatrix$ is computed based on the test stress following \refeq{eq:matrixstress}. 
This is used to compute the second deviatoric stress invariant $\JtwoMatrix$ for the average matrix stress.
Finally the effective compressive strength is estimated as
\begin{align}
	\fcEff = \frac{\JtwoTest}{\JtwoMatrix} \forceTest.
\end{align}
\subsubsection{Approximation of thermal conductivity}\label{ssec:thermalconductivity}
Homogenization of the thermal conductivity is based on the Mori-Tanaka method as well.
The formulation is similar to \refeq{eq:keff} and \refeq{eq:geff}.
The expressions are taken from \cite{str_2011_mbeo}.
The thermal conductivity $\thermCondHom$ is computed as
\begin{align}
	\thermCondHom &= \dfrac{\volFracMatrix\thermCondMatrix + \volFracIncl \thermCondIncl \concentrationThermCondIncl}{\volFracMatrix +  \volFracIncl \concentrationThermCondIncl}\quad\text{and}\\
	\concentrationThermCondIncl &= \frac{3\thermCondMatrix}{2\thermCondMatrix+\thermCondIncl}.
\end{align}

\section{FE Concrete Model}\label{appendix:fem}
\subsection{Modeling of the temperature field}
The temperature distribution is generally described by the heat equation as
\begin{align}
	\density \heatCapSpecific \dTdt = \nabla \cdot (\thermCondEff \nabla \temp) + \dQdt \label{eq:heat1}
\end{align}
with $\thermCondEff$ the effective thermal conductivity, $C$ the specific heat capacity, $\density$ the density and $\density \heatCapSpecific$ the volumetric heat capacity. The volumetric heat $Q$ due to hydration is also called the latent heat of hydration, or the heat source. In the paper, the density, the thermal conductivity and the volumetric heat capacity to be constants assumed to be sufficiently accurate for our purpose, even though there are more elaborate models taking into account effects of  temperature, moisture and/or the hydration.

\subsubsection{Degree of hydration\texorpdfstring{ $\DOH$}{}}
The degree of hydration $\DOH$ is defined as the ratio between the cumulative heat $\heat$ at time $\zeit$ and the total theoretical volumetric heat by complete hydration $\heatInf$:
\begin{align}
	\DOH(\zeit) = \frac{\heat(\zeit)}{\heatInf} \label{eq:doh}
\end{align}
assuming a linear relation between the degree of hydration and the heat development. Therefore, the time derivative of the heat source $\dot{\heat}$ can be rewritten in terms of $\DOH$, 
\begin{align}
	\dQdt = \dDOHdt \heatInf. \label{eq:qdotalphadot}
\end{align}
Approximated values for the total potential heat range between 300 and 600$\frac{J}{g}$ for binders of different cement types, e.g. Ordinary Portland cement $\heatInf= $ 375–525$\frac{J}{g}$ or Pozzolanic cement $\heatInf= $ 315–420$\frac{J}{g}$.  

\subsubsection{Affinity}
\begin{figure*}[bt]%
	\centering
	\includegraphics[width=\textwidth]{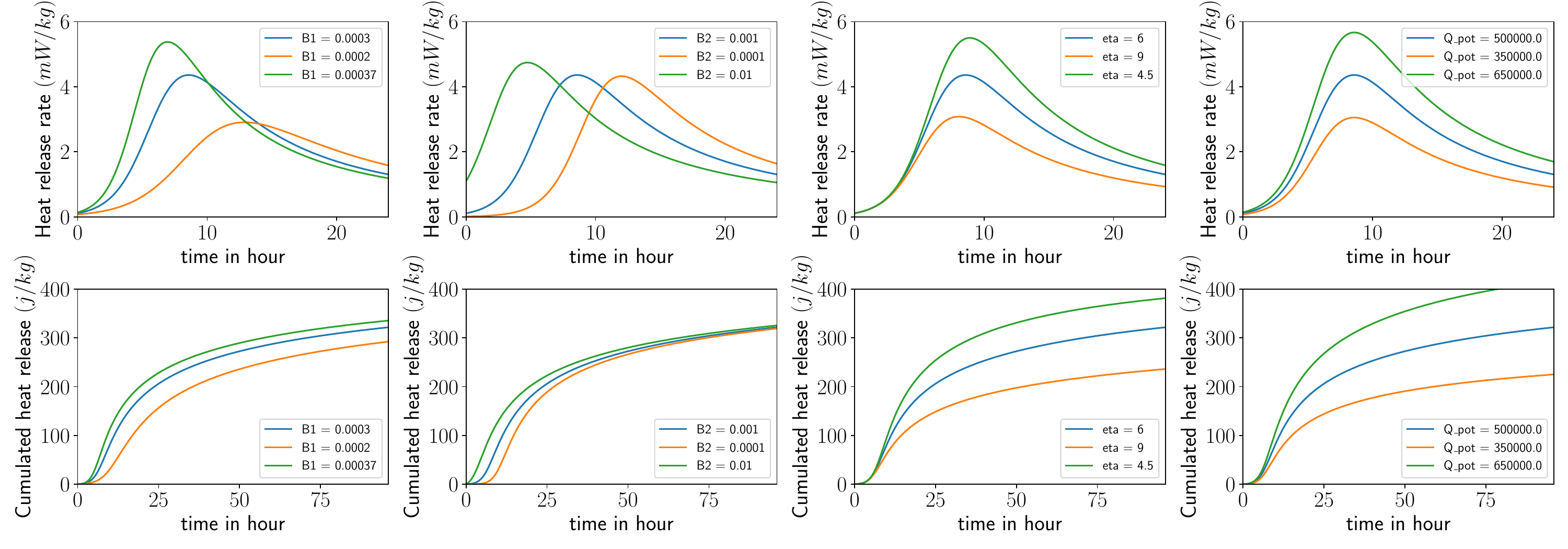}
	\caption{Influence of the hydration parameters on the heat release rate and the cumulative heat release.}\label{fig:heatrelease}
\end{figure*}

The heat release can be modeled based on the chemical affinity $\affinity$ of the binder. The hydration kinetics are defined as a function of affinity at a reference temperature $\affinityTemp$ and a temperature dependent scale factor ${\affinityScale}$
\begin{align}
	\dot{\DOH} = \affinityTemp(\DOH){\affinityScale}(\temp),\label{eq:affinitydot}
\end{align}
The reference affinity, based on the degree of hydration, is approximated by
\begin{align}
	\affinityTemp(\DOH) = \hydParBone\!\!\left(\frac{\hydParBtwo}{\DOHmax} + \DOH\right)\!\!(\DOHmax - \DOH)\exp\!\!\left(\!\!-\hydParEta \frac{\DOH}{\DOHmax}\right),
\end{align}
where $\hydParBone$ and $\hydParBtwo$ are coefficients depending on the binder.
The scale function is given as
\begin{align}
	\affinityScale = \exp\left(-\frac{\activE}{\gasConst}\left(\dfrac{1}{\temp}-\dfrac{1}{\tempRef}\right)\right).
\end{align}
An example function to approximate the maximum degree of hydration based on the water to cement mass ratio $\wc$, by Mills (1966)
\begin{align}
	\DOHmax = \dfrac{1.031\wc}{0.194 + \wc},
\end{align}
this refers to Portland cement. 
\reffig{fig:heatrelease} shows the influence of the three numerical parameters $B_1$, $B_2$, $\eta$ and the potential heat release $\heatInf$ on the heat release rate as well as on the cumulative heat release .

\subsubsection{Discretization and solution}
Using \refeq{eq:qdotalphadot} in \refeq{eq:heat1}, the heat equation is given as
\begin{align}
	\density \heatCapSpecific \frac{\partial \temp}{\partial \zeit} = \nabla \cdot (\thermCondEff \nabla \temp) + \heatInf\frac{\partial  \DOH}{\partial \zeit} 
\end{align}
Now we apply a backward Euler scheme 
\begin{align}
	\dot{\temp} =& \dfrac{\tempCurrent-\tempLast}{\Delta \zeit} \quad\text{and}\label{eq:timediscr}\\
	\dot{\DOH} =& \dfrac{\Delta \DOH}{\Delta \zeit}  \quad\text{with}\quad
	\Delta \DOH  = \DOHCurrent- \DOHLast
	\label{eq:timediscr2},
\end{align}
drop the index $\currentn$ for readability to obtain
\begin{align}
	\density \heatCapSpecific \temp  - {\Delta \zeit} \nabla \cdot (\thermCondEff \nabla \temp) - \heatInf\Delta \DOH
	= \density \heatCapSpecific \tempLast \label{eq:heat2}.
\end{align}
Using \refeq{eq:timediscr2} and \refeq{eq:affinitydot}, a formulation for $\Delta \DOH$ is obtained
\begin{align}
	\Delta \DOH = \Delta \zeit \affinityTemp( \DOH)\affinityScale(\temp). \label{eq:deltaalpha}
\end{align}
We define the affinity in terms of $ \DOHLast$ and $\Delta \DOH$ to solve for $\Delta \DOH$ on the quadrature point level 
\begin{align}
	\affinityTemp =& \hydParBone\exp\left(-\hydParEta \frac{\Delta \DOH+\DOHLast}{\DOHmax}\right)
	\left(\tfrac{\hydParBtwo}{\DOHmax} + \Delta\DOH+\DOHLast\right)\cdot\nonumber\\ &(\DOHmax - \Delta\DOH - \DOHLast).
\end{align}
Now we can solve the nonlinear function 
\begin{align}
	\function(\Delta\DOH) = \Delta\DOH - \Delta \zeit \affinityTemp(\Delta\DOH)\affinityScale(\temp) = 0
\end{align}
using an iterative Newton-Raphson solver.

\subsection{Coupling Material Properties to Degree of Hydration} \label{appendix:fem_evolution}
\subsubsection{Compressive strength}
The compressive strength in terms of the degree of hydration can be approximated using an exponential function, c.f. \cite{car_2016_mamt},
\begin{align}
	\fc(\DOH) = \DOH(\zeit)^{\strengthCExp} \fcInf \label{eq:fcwrtDOH}.
\end{align}
This model has two parameters, $\fcInf$, the compressive strength of the parameter at full hydration, $\DOH = 1$ and $\strengthCExp$ the exponent, which is a material parameter that characterizes the temporal evolution.

The first parameter could theoretically be obtained through experiments.
However the total hydration can take years. Therefore, we can compute it using the 28 days values of the compressive strength and the corresponding degree of hydration
\begin{align}
\fcInf = 	\dfrac{\fcTwentyEight}{{\DOHTwentyEight}^{\strengthCExp}}.
\end{align}
\subsubsection{Young's Modulus}
The publication \cite{car_2016_mamt} proposes a model for the evolution of the Young's modulus assuming an initial linear increase of the Young's modulus up to a degree of hydration $\DOHt$,
\begin{align}
	\eMod(\DOH) = 
	\begin{cases}
		\eModInf \frac{\DOH(\zeit)}{\DOHt}{\DOHt}^{\stiffExp}   
		& \text{for $\DOH < \DOHt$}\\
		\eModInf {\DOH(\zeit)}^{\stiffExp}  
		& \text{for $\DOH \ge \DOHt$}.
	\end{cases}\label{eq:EwrtDOH}
\end{align}
Contrary to other publications, no dormant period is assumed.
Similarly to the strength standardized testing of the Young's modulus is done after 28 day, $\eModTwentyEight$.
To effectively use these experimental values, $\eModInf$ is approximated as
\begin{align}
	\eModInf = \dfrac{\eModTwentyEight}{{\DOHTwentyEight}^{\stiffExp}}.
\end{align}
using the approximated degree of hydration.

\subsection{Constraints}\label{sec:appendix_constraints}
The FEM simulation is used to compute two practical constraints relevant to the precast concrete industry.
At each time step, the worst point is chosen to represent the part, therefore ensuring that the criterion is fulfilled in the whole domain.
The first constraint limits the maximum allowed temperature.
The constraint is computed as the normalized difference between the maximum temperature reached $\tempMax$ and the temperature limit $\tempLimit$ 
\begin{align}
	\FEMConstraintT = \frac{\tempMax - \tempLimit}{\tempLimit}, \label{eq:concstraintT}
\end{align}
where $\FEMConstraintT > 0$ is not admissible, as the temperature limit \inputtemperaturelimit \degree C has been exceeded.\\

The second constraint is the estimated time of demolding.
This is critical, as the manufacturer has a limited number of forms.
The faster the part can be demolded, the faster it can be reused, increasing the output capacity.
On the other hand, the part must not be demolded too early, as it might get damaged while being moved.
To approximate the minimal time of demolding, a constraint is formulated based on the local stresses $\FEMConstraintStress$.
It evaluates the Rankine criterion for the principal tensile stresses, using the yield strength of steel $\beamfs$ and a simplified Drucker-Prager criterion, based on the evolving compressive strength of the concrete $\fc$,
\begin{align}
	\FEMConstraintStress =  \max
	\begin{cases}
		\FEMConstraintRK = \frac{\|\principalStressTension\| - \beamfs}{\beamfs} \\
		\FEMConstraintDP =  \frac{\sqrt{\frac{1}{3} \firstStressInvariant^2 - \secondStressInvariant} - \frac{\fc^3}{\sqrt{3}}}{\fc}
	\end{cases},\label{eq:constraintStress}
\end{align}
where $\FEMConstraintStress > 0$ is not admissible.
In contrast to standard yield surfaces, the value is normalized, to be unit less.
This constraint aims to approximate the compressive failure often simulated with plasticity and the tensile effect of reinforcement steel.
As boundary conditions, a simply supported beam under it own weight has been chosen to approximate possible loading condition while the part is moved.
This constraint is evaluated for each time step in the simulation.
The critical point in time is approximated where $\FEMConstraintStress(t_{crit}) = 0$. This is normalized with the prescribed time of demoulding to obtain a dimensionless constraint.

\section{Beam Design}\label{appendix:beam}
Following design code \citeauthor{DIN1992-1-1} for a singly reinforced beam, meaning a reinforced concrete beam with only reinforcement at the bottom.
The assumed cross section is rectangular

\subsection{Maximum bending moment}
Assuming a simply supported beam with a given length $\beamLength$ in mm, a distributed load $\beamDistrLoad$ in N/mm and a point load $\beamPointLoad$ in N/mm
the maximum bending moment $\beamMaxMoment$ in N/mm$^2$ is computed as
\begin{align}
	\beamMaxMoment= \beamDistrLoad \frac{\beamLength^2}{8} + \beamPointLoad \frac{\beamLength}{4}.
\end{align}
The applied loads already incorporate any required safety factors.
\subsection{Computing the minimal required steel reinforcement}
Given a beam with the height $\beamHeight$ in mm, a concrete cover of $\beamCover$ in mm, a steel reinforcement diameter of $\beamSteelDiameter$ in mm for the longitudinal bars and a bar diameter of $\beamSteelDiameter$ in mm for the transversal reinforcement also called stirrups,
\begin{align}
	\beamHeightEff = \beamHeight - \beamCover - \beamSteelDiameterStirrups - \frac{1}{2} \beamSteelDiameter.
\end{align}
According to the German norm standard safety factors are applied, $\beamTimeSF = 0.85$, $\beamConcreteSF = 1.5$ and $\beamSteelSF = 1.15$, leading to the design compressive strength for concrete $\beamfcd$ and  the design tensile yield strength $\beamfsd$ for steel
\begin{align}
	\beamfcd = \beamTimeSF \frac{\fc }{\beamConcreteSF}\\
	\beamfsd = \frac{\beamfs}{\beamSteelSF},
\end{align}
where $\fc$ denotes the concrete compressive strength and $\beamfs$ the steel's tensile yield strength.\\
To compute the force applied in the compression zone, the lever arm of the applied moment $\beamLeverMoment$ is given by 
\begin{align}
	\beamLeverMoment &= \beamHeightEff(0.5+\sqrt{0.25-0.5 \beamK}),\quad \text{with}
	\\
	\beamK &= \frac{\beamMaxMoment }{\beamWidth \beamHeightEff^2 \beamfcd}. \label{eq:compr_str_constr}
\end{align}
The minimum required steel $\beamSteelReq$ is then computed based on the lever arm, the design yield strength of steel and the maximum bending moment, as
\begin{align}
	\beamSteelReq = \frac{\beamMaxMoment}{\beamfsd \beamLeverMoment}.\label{eq:Areq}
\end{align}
\subsection{Optimization constraints}
\subsubsection{Compressive strength constraint}
Based on \refeq{eq:compr_str_constr}, we define the compressive strength constraint as
\begin{align}
	\beamConstraintFc = \beamK - 0.5, \label{eq:constraintfc}
\end{align}
where $\beamConstraintFc > 0$ is not admissible, as there is no solution for \refeq{eq:compr_str_constr}.
\subsubsection{Geometrical constraint}
The geometrical constraint checks that the required steel area $\beamSteelReq$ does not exceed the maximum steel area $\beamSteelMax$ that fits inside the design space.
For our example, we assume the steel reinforcement is only arranged in a single layer.
This limits the available space for rebars in two ways, by the required minimal spacing $\beamMinSpacing$ between the bars, to allow concrete to pass, and by the required space on the outside, the concrete cover $\beamCover$ and stirrups diameter $\beamSteelDiameterStirrups$.
To compute $\beamSteelMax$, the maximum number for steel bars $\beamNSteelMax$ and the maximum diameter $\beamSteelDiameterMax$ from a given list of admissible diameters are determined that fulfill
\begin{align}
	s \ge& \beamMinSpacing,\quad\text{with} \\
	s =& \frac{\beamWidth - 2 \beamCover - 2 \beamSteelDiameterStirrups - \beamNSteelMax \beamSteelDiameterMax}{\beamNSteelMax - 1} \quad\text{and}\\
	\beamNSteelMax \ge& 2.
\end{align}
According to \citeauthor{DIN1992-1-1}, the minimum spacing between two bars $\beamMinSpacing$ is given by the minimum of the concrete cover (\inputconcretecover \,\text{\inputconcretecoverunit}) and the rebar diameter.
The maximum possible reinforcement is given by
\begin{equation}
	\beamSteelMax = \beamNSteel \pi \left( \frac{\beamSteelDiameter}{2}\right)^{2}.
\end{equation}
The geometry constraint is computed as
\begin{equation}
	\beamConstraintGeometry = \frac{\beamSteelReq - \beamSteelMax}{\beamSteelMax} \label{eq:constraintGeo}
\end{equation}
where $\beamConstraintGeometry > 0$ is not admissible, as the required steel area exceeds the available space.
\subsubsection{Combined beam constraint}
To simplify the optimization procedure, the two constraints are combined into a single one by using the maximum value,
\begin{equation}\label{eq:design_cons}
	\beamConstraintBeam = \max ( \beamConstraintGeometry, \beamConstraintFc).
\end{equation}
Evidently, this constraint is also defined as: $\beamConstraintBeam > 0$ is not admissible.

\section{Parameter Tables}\label{appendix:parameters}
This is the collection of the used parameter for the various example calculation.

\begin{table}[!htpb]
	\begin{center}
		\begin{minipage}{.5\textwidth}
			\caption{Parameters of the simply supported beam for the computation of the steel reinforcement}\label{tab:beamdesigninput}
			\begin{tabular}{lrl}
				\toprule
				Name &  Value&Unit\\
				\midrule
				Length & 1000 &$\beamExSpanUnit$\\
				Width & 350 &$\beamExWidthUnit$\\
				Height& \beamExHeightC &$\beamExHeightCUnit$\\
				Steel yield strength& 300 &$\beamExYieldStrSteelUnit$\\
				Diameter stirrups& \beamExSteelDiaBu &$\beamExSteelDiaBuUnit$\\
				Minimal concrete cover& \beamExCoverMin &$\beamExCoverMinUnit$\\
				Load& \beamExPointLoadC &$\beamExPointLoadCUnit$\\
				Concrete compressive strength& \beamExComprStrConcreteC &$\beamExComprStrConcreteCUnit$\\
				\botrule
			\end{tabular}
		\end{minipage}
	\end{center}
\end{table}

\end{appendices}
\end{document}